\documentclass[12pt]{amsart} 
\usepackage{amsmath,amssymb,amsthm,  mathtools} 
\usepackage[dvipsnames]{xcolor}
\usepackage{color,colordvi,graphicx,tikz}
\usepackage{ulem}  

\newtheorem{theorem}{Theorem}

\newtheorem{proposition}[theorem]{Proposition} 
\newtheorem{conjecture}[theorem]{Conjecture}

\theoremstyle{definition} 
\newtheorem{definition}[theorem]{Definition} 
\newtheorem{example}[theorem]{Example} 
\newtheorem{remark}[theorem]{Remark} 
\newcommand{\bAA}{{\mathbb A}}
\newcommand{\CC}{{\mathbb C}}
\newcommand{\EE}{{\mathbb E}}
\newcommand{\FF}{{\mathbb F}}
\newcommand{\GG}{{\mathbb G}}
\newcommand{\KK}{{\mathbb K}}
\newcommand{\LL}{{\mathbb L}}
\newcommand{\MM}{{\mathbb M}}

\newcommand{\PP}{{\mathbb P}}
\newcommand{\QQ}{{\mathbb Q}}
\newcommand{\RR}{{\mathbb R}}

\newcommand{\ZZ}{{\mathbb Z}}

\newcommand{\calA}{{\mathcal A}}

\newcommand{\calH}{{\mathcal H}}
\newcommand{\calO}{{\mathcal O}}
\newcommand{\calT}{{\mathcal T}}
\newcommand{\calV}{{\mathcal V}}

\newcommand{\bbk}{{\Bbbk}}

\newcommand{\bfI}{\mathbf{I}}
\newcommand{\bfk}{\mathbf{k}}
\newcommand{\bfp}{\mathbf{p}}
\newcommand{\bfR}{\mathbf{R}}
\newcommand{\bft}{\mathbf{t}}
\newcommand{\bfv}{\mathbf{v}}
\newcommand{\bfx}{\mathbf{x}}
\newcommand{\bfy}{\mathbf{y}}
\newcommand{\bfzero}{\mathbf{0}}
\newcommand{\SE}{\text{SE}}
\newcommand{\SO}{\text{SO}}

\newcommand{\conv}{\mbox{\rm conv}}
\newcommand{\rank}{\mbox{\rm rank}}
\newcommand{\MV}{\mbox{\rm MV}}

\newcommand{\Gal}{\mbox{\rm Gal}}
\newcommand{\Hom}{\mbox{\rm Hom}}
\newcommand{\Mon}{\mbox{\rm Mon}}
\newcommand{\Gr}{\mbox{\rm Gr}}

\newcommand{\Fdot}{{F^\bullet}}
\newcommand{\adot}{{a_\bullet}}
\newcommand{\calAdot}{{\mathcal A}_\bullet}
\newcommand{\kapDot}{{\kappa^\bullet}}
\newcommand{\lamDot}{{\lambda^\bullet}}
\newcommand{\muDot}{{\mu^\bullet}}

\newcommand{\Fl}{{\mathbb F}{\ell}}

\newcommand{\cbd}{{\color{black}\bullet}{\color{red}\bullet}}

\newcommand{\defcolor}[1]{{\color{blue}#1}} 
\newcommand{\demph}[1]{\defcolor{{\sl #1}}} 

\newcounter{FNC}[page]
\def\fauxfootnote#1{{\addtocounter{FNC}{2}${\color{magenta}^\fnsymbol{FNC}}$%
     \let\thefootnote\relax\footnotetext{{\color{magenta}$^\fnsymbol{FNC}$#1}}}}

\title[Galois Groups in Enumerative Geometry and Applications]{Galois Groups in Enumerative Geometry\\ and Applications} 
\author[F.~Sottile]{Frank Sottile} 
\address{F.~Sottile\\ 
         Department of Mathematics\\ 
         Texas A\&M University\\ 
         College Station\\ 
         Texas \ 77843\\ 
         USA} 
\email{sottile@tamu.edu} 
\urladdr{https://franksottile.github.io/} 
\author[T.~Yahl]{Thomas Yahl} 
\address{T.~Yahl\\ 
         Department of Mathematics\\ 
         University od Wisconsin\\
         Madison\\
         Wisconsin \ 53706\\ 
         USA} 
\email{tyahl@wisc.edu} 
\urladdr{https://tjyahl.github.io/}
\thanks{Research of Sottile and Yahl supported by grant 636314 from the Simons Foundation,
        and the National Science Foundation through grant  DMS-2201005.}
\subjclass{11R32, 12F12, 14N15, 14M25, 14G99, 14Q65, 65H14}
\keywords{Galois group, enumerative geometry, sparse polynomial systems, Schubert caculus, Fano problem, homotopy continuation} 
\begin{document} 
\begin{abstract}
  As Jordan observed in 1870, just as univariate polynomials have Galois groups, so do problems in enumerative geometry.
  Despite this pedigree, the study of Galois groups in enumerative geometry was dormant for a century, with a systematic study only
  occurring in the past 15 years.
  We discuss the current directions of this study, including open problems and conjectures.
\end{abstract}
\maketitle

%
\section*{Introduction}

We are all familiar with Galois groups:  They play an important role in the structure of field extensions
and control the solvability of equations.
Less known is that they have a long history in enumerative geometry.
In fact, the first comprehensive treatise on Galois theory, Jordan's ``Trait\'e des Substitutions et des \'Equations
alg\'ebriques''~\cite[Ch.~III]{J1870}, also discusses Galois theory in the context of several classical problems in enumerative geometry.

While Galois theory developed into a cornerstone of number theory and of arithmetic geometry, its role in enumerative geometry lay 
dormant until Harris's 1979 paper ``Galois groups of enumerative problems''~\cite{Harris}.
Harris revisited Jordan's treatment of classical problems and gave a proof that, over $\CC$, the Galois and monodromy groups coincide.
He used this to introduce a geometric method to show that an enumerative Galois group is the full symmetric group and showed that
several enumerative Galois groups are full-symmetric, including  generalizations of the classical problems studied by Jordan.

We sketch the development of Galois groups in enumerative geometry since 1979.
This includes some new and newly applied methods to study or compute Galois groups in this context, as well as recent results 
and open problems.
A theme that Jordan initiated is that intrinsic structure of the solutions to an enumerative problem constrains its Galois group
\defcolor{$G$} giving an ``upper bound'' for $G$.
The problem of identifying the Galois group $G$ becomes that of showing it is as ``large as possible''.
In all cases when $G$ has been determined, it is as large as possible given the intrinsic structure.
Thus we may view $G$ as encoding the
intrinsic structure of the enumerative problem.

\begin{figure}[htb]
  \centering
  \begin{picture}(110,120)
       \put(0,0){\includegraphics[height=130pt]{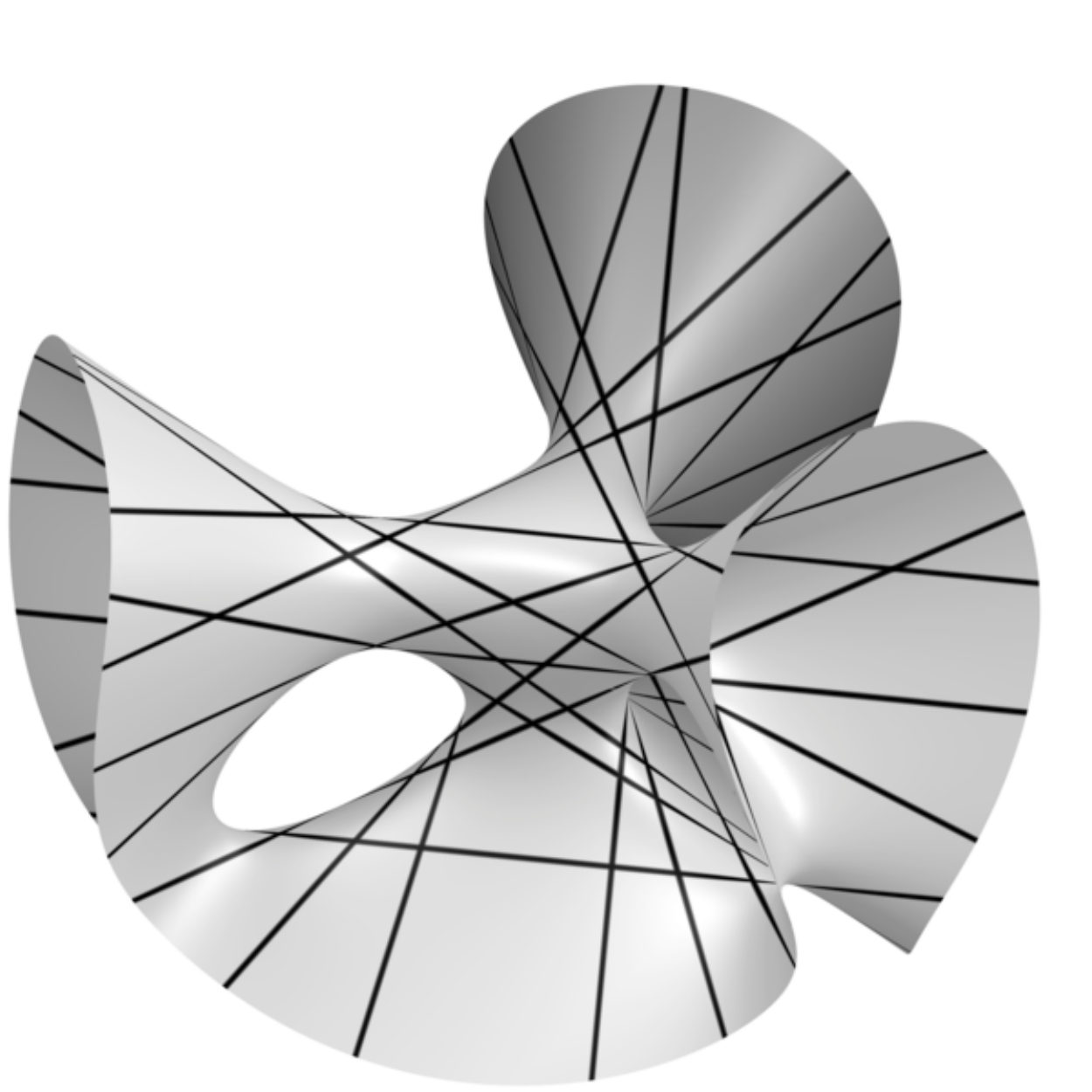}}
  \end{picture}
  \caption{A cubic with 27 lines. (Image courtesy of Oliver Labs)}
  \label{F:cubic}
\end{figure}

Consider the problem of lines on a cubic surface.
Cayley~\cite{Cayley1849} and Salmon~\cite{Salmon1849} showed that a smooth cubic surface $\calV(f)$ in $\PP^3$ ($f$ is a homogeneous
cubic in four variables) contains 27 lines.
(See Figure~\ref{F:cubic}.)
 This holds over any algebraically closed field.
%
%
%
When $f$ has rational coefficients, the field \defcolor{$\KK$} of definition of the lines is a Galois extension of $\QQ$, and its Galois group
\defcolor{$G$} has a faithful action on the 27 lines.

As the lines lie on a surface, we expect that some will meet, and Schl\"afli~\cite{Schlafli1858} showed that for a general
cubic, these lines form a remarkable incidence configuration whose symmetry group is
what we now understand to be the Coxeter
group $E_6$.
%
%
As Jordan observed, this implies that $G$ is a subgroup of $E_6$, and it is now known that for most cubic surfaces $G=E_6$.

A modern view begins with the incidence variety of this enumerative problem.
The space of homogeneous cubics on $\PP^3$ forms a 19-dimensional projective space, as a cubic in four variables has $\binom{3+4-1}{3}=20$
coefficients. 
Writing \defcolor{$\GG(1,\PP^3)$} for the (four-dimensional) Grassmannian of lines in $\PP^3$, we have the incidence
variety, 
 \begin{equation}\label{Eq:Cubic_Family}
   \raisebox{-30pt}{\begin{picture}(232,65)(-2,-2)
   \put(0,52){$\defcolor{\Gamma}\ 
         \vcentcolon=\ \{(\ell,f)\in \GG(1,\PP^3)\times\PP^{19}_{\mbox{\tiny cubics}}\,\mid\, f|_\ell\equiv 0\}$\,.}
    \put(2,48){\vector(0,-1){32}}  \put(6,32){\small$\pi$}
    \put(-1,2){$\PP^{19}_{\mbox{\tiny cubics}}$}
   \end{picture}}
 \end{equation}
Let \defcolor{$\bbk$} be our ground field, which we assume for now to be algebraically closed.
Both $\Gamma$ and $\PP^{19}$ are irreducible.
Consider their fields of rational functions, $\bbk(\Gamma)$ and $\bbk(\PP^{19})$.
As the typical fiber of $\pi$ consists of 27 points and $\pi$ is dominant, $\pi^*(\bbk(\PP^{19}))$ is a subfield of $\bbk(\Gamma)$, and the
extension has degree 27.
The Galois group \defcolor{$G$} of the normal closure of this extension acts on the lines in the generic cubic surface over $\PP^{19}$, and 
we have that $G=E_6$.

Suppose now that $\bbk=\CC$.
If $B\subset\PP^{19}$ is the set of singular cubics (a degree 32 hypersurface~\cite{Salmon1860}) then over $\PP^{19}{\smallsetminus}B$,
$\Gamma$ is a covering space of degree 27.
Lifting based loops gives the monodromy action of the fundamental group of $\PP^{19}{\smallsetminus}B$ on the fiber above the base point.
Permutations of the fiber obtained in this way constitute the \demph{monodromy group} of $\pi$.
For the same  reasons as before, this is a subgroup of $E_6$.
In fact, it equals $E_6$.

This situation, a dominant map $\pi\colon X\to Z$ of irreducible equidimensional varieties, is called a
\demph{branched  cover}. 
Branched covers are common in enumerative geometry and applications of algebraic geometry.
For the problem of 27 lines, that the algebraic Galois group equals the geometric monodromy group is no accident;
While Harris~\cite{Harris} gave a modern proof, the equality of these two groups may be traced back to
Hermite~\cite{Hermite}. 
We sketch a proof, valid over arbitrary fields, in Section~\ref{Sec:background}.

Harris's article brought this topic into contemporary algebraic geometry.
He also introduced geometric methods to show that the Galois group of an enumerative problem is
\demph{fully  symmetric} in that it is the full symmetric group on the solutions.
In the 25 years following its publication, the Galois group
was determined in only a handful of enumerative problems.
For example, D'Souza~\cite{DS88} showed that the problem of lines in $\PP^3$ tangent to a smooth octic surface at four points
(everywhere tangent lines) had Galois group that is fully symmetric. 
Interestingly, he did not determine the number of everywhere tangent lines.

This changed in 2006 when Vakil introduced a method~\cite{Va06b} to deduce that the Galois group of a Schubert problem on a
Grassmannian (a \demph{Schubert Galois group}) contains the alternating group on its solutions.
Such a Galois group is
said to be \demph{giant}, as other proper subgroups of a symmetric group are relatively minuscule.
Vakil used his method
to show that most Schubert problems on small Grassmannians were 
giant, and to discover an
infinite family of Schubert problems whose Galois groups were not the full symmetric group.
As we saw in the problem of 27 lines on a cubic surface, such an enumerative problem with a small Galois group typically possesses
  some internal structure.
Consequently, we use the adjective \demph{enriched} to describe such a problem or Galois group.
Enriched Schubert problems were also found on more general flag manifolds~\cite{RSSS}.
These discoveries inspired a more systematic study of Schubert Galois groups, which we discuss in Section~\ref{S:SchubertGaloisGroups}.
Despite significant progress, the \demph{inverse Galois problem} for Schubert calculus remains open.

Galois groups of enumerative problems are usually transitive permutation groups.
There is a dichotomy between those transitive permutation groups that preserve no nontrivial partition, called \demph{primitive} groups, and
the \demph{imprimitive} groups that do preserve a nontrivial partition. 
The Galois group of the 27 lines is primitive, but most known enriched Schubert problems have imprimitive Galois groups.

Another well-understood class of enumerative problems comes from the Bernstein-Kushnirenko Theorem~\cite{Bernstein75,Kushnirenko}.
This gives the number of solutions to a system of polynomial equations that are general given the monomials
occurring in the equations. 
Esterov~\cite{Esterov} determined which of these problems have fully symmetric Galois group and showed that all others have an imprimitive
Galois group.
Here, too, the inverse Galois problem remains open.
We discuss this in Section~\ref{S:GGsparse}.

The problem of lines on a cubic surface is the first in the class of \demph{Fano problems}, which involve counting the number of linear
subspaces that lie on a general complete intersection in projective space.
Recently, Hashimoto and Kadets~\cite{HK} nearly determined the Galois groups of all Fano problems.
Most are 
giant, except for the lines on a cubic surface and the $r$-planes lying on the intersection of two quadrics in
$\PP^{2r+2}$.
We explain this in Section~\ref{Sec:Fano}, and discuss computations from~\cite{Yahl} which show that several Fano
problems of moderate size are full-symmetric. 

Branched covers arise from families of polynomial systems, which are common in the applications of mathematics.
Oftentimes the application or the formulation as a system of polynomials possesses some intrinsic structure, which is manifested in the
corresponding Galois group being enriched.
In Section~\ref{Sec:Applications}, we discuss two occurrences of enriched Galois groups in applications and a
  computational method that exploits structure in Galois groups for computing solutions to systems of equations.

We begin in Section~\ref{Sec:background} with a general discussion of Galois groups in enumerative geometry, and sketch some methods
  from numerical algebraic geometry in Section~\ref{Sec:NAG}.
Later, in Section~\ref{Sec:Computing}, we present methods, both numerical and symbolic, to compute and study Galois groups in enumerative
geometry. 

%
\section{Galois groups of branched covers}
\label{Sec:background}

We will let \defcolor{$\bbk$} be a field with algebraic closure \defcolor{$\overline{\bbk}$}, which we fix.
We adopt standard terminology from algebraic geometry:
An affine (projective) scheme \defcolor{$\calV(F)$} is defined in $\bAA^n$ ($\PP^n$) by polynomials (homogeneous forms)
$F=(f_1,\dotsc,f_m)$ in $n$ ($n{+}1$) variables with coefficients in $\bbk$.
We will call the collection $F$ a \defcolor{system} (of equations) and say the isolated points of $\calV(F)$
(over $\overline{\bbk}$) are the \demph{solutions} to $F$. 
The scheme $\calV(F)$ is a \demph{variety} when every irreducible component of $\calV(F)$ is reduced.
We may also use variety to refer to the underlying variety.
We write \defcolor{$X(\overline{\bbk})$} for the points of a variety $X$ with coordinates in $\overline{\bbk}$.

Recall that the Galois group of a separable univariate polynomial $f(x)\in \bbk[x]$ is the Galois group of the splitting
field of $f$, which is generated over $\bbk$ by the roots of $f$ and is a subfield of $\overline{\bbk}$.
Given a system $F$ of multivariate polynomials over $\bbk$, its splitting field is the
subfield of $\overline{\bbk}$ generated over $\bbk$ by the coordinates of all
solutions to $F$, and its Galois group is the Galois group of this field extension.

A separable map $\pi\colon X\to Z$ of irreducible varieties is a \demph{branched cover} when $X$ and $Z$ have the same dimension and
$\pi(X)$ is dense in $Z$ ($\pi$ is \demph{dominant}).
Branched covers are ubiquitous in enumerative geometry and in applications of algebraic geometry.
When the varieties are complex, there is a proper subvariety $B\subset Z$ (the \demph{branch locus}) such that $\pi$ is a covering space
over $Z{\smallsetminus} B$.
We explain how to associate a Galois/monodromy group to a branched cover and then give some background on permutation groups, and the relation
between imprimitivity of the Galois group and decomposability of the branched cover.

\subsection{Galois and monodromy groups of branched covers}\label{SS:Galois_Monodromy}

Let $\pi\colon X\to Z$ be a branched cover.
As $\pi$ is dominant, the function field \defcolor{$\bbk(Z)$} of $Z$  embeds as a subfield of the function field $\bbk(X)$ of $X$.
This realizes $\bbk(X)$ as a finite extension of $\bbk(Z)$ of degree $d$, the \demph{degree of $\pi$}.
Let \defcolor{$\KK$} be the normal closure of this extension.
The \defcolor{Galois group} of the branched cover $\pi$, denoted \defcolor{$\Gal_\pi$}, is the Galois group of $\KK/\bbk(Z)$.
This is a transitive subgroup of the symmetric group $S_d$ that is well-defined up to conjugation.

There is also a geometric construction of $\Gal_\pi$.
For $1\leq s\leq d$, let \defcolor{$X_Z^s$} be the $s$-th fold iterated fiber product of $\pi\colon X\to Z$,
  \[
    X_Z^s\ \vcentcolon=\ \underbrace{X\times_Z X\times_Z \dotsb \times_Z X}_{s}\,.
  \]
(This is the pull-back of the $s$-fold Cartesian product $X^s\to Z^s$  along the diagonal map
    $\Delta\colon Z\hookrightarrow Z^s$.)
The fiber of $\pi^s\colon X_Z^s\to Z$ over a point $z\in Z$ is the $s$-fold Cartesian product $(\pi^{-1}(z))^s$ of the
fiber of $\pi$ over $z$. 

The fiber product has many irreducible components when $s>1$, possibly of different dimensions.
Let $\defcolor{U}\subset Z$ be the maximal dense open subset over which $\pi$ is proper and
\'etale---fibers $\pi^{-1}(z)$ for $z\in
U$ are zero-dimensional reduced schemes of degree $d$.
Its complement is the \demph{branch locus} \defcolor{$B$} of $\pi$.
The \demph{big diagonal} of $X_Z^s$ is the closed subscheme consisting of $s$-tuples with a repeated coordinate.
Let $\defcolor{X_Z^{(s)}}$ be the closure in $X_Z^s$ of the complement of the big diagonal in $(\pi^s)^{-1}(U)$.
The fiber of $X_Z^{(s)}$ over a point $z\in U(\overline{\bbk})$ consists of $s$-tuples of distinct points of the fiber $\pi^{-1}(z)$.

When $s=d$, the symmetric group $S_d$ acts on $X_Z^{(d)}$, permuting each $d$-tuple.
It permutes the irreducible components and acts simply transitively on the fiber above a point $z\in U(\overline{\bbk})$.
Let $X'\subset X_Z^{(d)}$ be an irreducible component (they are all isomorphic when $s=d$).

We compare this  to the construction of the splitting field of a univariate polynomial.
Replacing $X$ and $Z$ by appropriate affine open subsets, we may embed $X$ as a hypersurface in  $Z\times \bAA^1_t$
with $X\to Z$ the projection.
Writing \defcolor{$\bbk[X]$} and \defcolor{$\bbk[Z]$} for their coordinate rings, there is a monic irreducible polynomial
$\defcolor{f}\in \bbk[Z][t]$ of degree $d$ such that $\bbk[X]=\bbk[Z][t]/\langle f\rangle$.
Thus $\bbk(X)=\bbk(Z)[t]/\langle f\rangle = \bbk(Z)(\alpha)$, where $\alpha$ is the image of $t$ in $\bbk[X]$.
If $X'$ is an irreducible component of $X_Z^{(d)}$, then $\bbk(X')=\bbk(Z)(\alpha_1,\dotsc,\alpha_d)$ where
$\alpha_i\in \bbk[X']$ is given by the composition of inclusion $X'\subset X^{(d)}_Z$, the $i$th coordinate projection $X^{(d)}_Z\to X$, and
the function $\alpha$.
As $i\neq j \Rightarrow \alpha_i\neq\alpha_j$ ($X'$ does not lie in the big diagonal), we see that $\alpha_1,\dotsc,\alpha_d$ are the roots
of $f$ in $\bbk(X')$.
Thus $\bbk(X')$ is the splitting field of $f$ and is Galois over $\bbk(Z)$.

The \demph{monodromy group} \defcolor{$\Mon_\pi$} of the branched cover is the subgroup of $S_d$ that preserves $X'$ (sends
points of $X'$ to points of $X'$).
Elements of $\Mon_\pi$ induce automorphisms of the extension $\bbk(X')/\bbk(Z)$ so that
$\Mon_\pi\subset\Gal(\bbk(X')/\bbk(Z))$, the Galois group of $\bbk(X')/\bbk(Z)$.
Since  $\Mon_\pi$ acts simply transitively on fibers of $X'\to Z$ above points in $U(\overline{\bbk})$, its order is the
degree of the map $X'\to Z$, which is the order of the field extension $\bbk(X')/\bbk(Z)$.
Hence we arrive at the result $\Mon_\pi=\Gal(\bbk(X')/\bbk(Z))$.

\begin{theorem}[Galois equals monodromy]
  \label{Th:Gal=Mon}
  For a branched cover $\pi\colon X\to Z$ defined over a field $\bbk$, the Galois group is equal to the monodromy group,
  \[
    \Gal_\pi\ =\ \Mon_\pi\,.
  \]
\end{theorem}

The enumerative problem of 27 lines on a cubic surface has a corresponding incidence
variety~\eqref{Eq:Cubic_Family} which is a branched cover, and its Galois/monodromy group is a special case of the results
of this section.
Incidence varieties of enumerative problems typically are branched covers and therefore have Galois groups as we will see
throughout this survey. 

We make an important observation.
While the Galois group of a branched cover $\pi\colon X\to Z$ is defined via a geometric construction, it does depend upon
the field of definition.
For example, consider the branched cover $\pi\colon\bAA^1\to\bAA^1$ given by $x\mapsto x^3$. 
Assume that $\bbk$ does not have characteristic 3, for otherwise $\pi$ is inseparable.
Over the rational numbers, this is $\pi\colon\bAA^1(\QQ)\to\bAA^1(\QQ)$, which has Galois group $S_3$.
Over any field containing $\sqrt{-3}$
(e.g.\ $-3$ is a square in $\bbk$) its Galois group is $A_3=\ZZ/3\ZZ$.
This is because the discriminant of the cubic $x^3-t$ defining $\pi$ is $-27t^2$, which is a square in
$\bbk(t)$ when $\sqrt{-3}\in\bbk$. 
When necessary, we write \defcolor{$\Gal_\pi(\bbk)$} to indicate that the branched cover is defined over $\bbk$.
While this notation uses the base field $\bbk$, it is important to keep in mind that this is a Galois group over the
transcendental extension $\bbk(Z)$ of $\bbk$.

We record some facts about how $\Gal_\pi(\bbk)$ behaves under field extensions.

\begin{theorem} \label{Th:GaloisUnderExtension}
  Suppose that $\pi\colon X\to Z$ is a branched cover defined over $\bbk$ and $\FF/\bbk$ is any field extension.
  Then $\Gal_\pi(\FF)$ is isomorphic to a subgroup of $\Gal_\pi(\bbk)$.

  When $\FF\supset \overline{\bbk}$ and $Z$ is a rational variety,  $\Gal_\pi(\FF)$ is isomorphic to a normal subgroup
  of $\Gal_\pi(\bbk)$. 
\end{theorem}
\begin{proof}
  Let $\KK/\bbk(Z)$ be the normal closure of the extension $\bbk(X)/\bbk(Z)$ and $\MM/\FF(Z)$ be the normal closure
  of $\FF(X)/\FF(Z)$.
Setting $\defcolor{\EE}\vcentcolon=\KK\cap\FF(Z)$, we have the following diagram of field extensions:
\begin{equation}\label{Eq:fieldExtensions}
   \raisebox{-38pt}{\begin{picture}(70,88)(2,-4)
            \put(33,78){$\MM$}
   \put( 9,59){\line(3,2){24}}     \put(55,65){\line(-1,1){10}} 
   \put( 2,47){$\KK$}               \put(49,54){$\FF(Z)$}
    \put(20,34){\line(-1,1){10}}    \put(30,34){\line(3,2){24}}
    \put(21,23){$\EE$}
    \put(25,11){\line(0,1){10}}
    \put(15, 0){$\bbk(Z)$}
    \put(76,38){.}
  \end{picture}}
\end{equation}
As $\MM$ is the compositum of $\KK$ and $\FF(Z)$,~\cite[Thm.\ VI.1.12]{Lang} implies that
 \begin{equation}\label{Eq:equalityOfGroups}
   \Gal_\pi(\FF)\ =\ \Gal(\MM/\FF(Z))\ =\ \Gal(\KK/\EE)\,.
 \end{equation}
Because $\EE$ is an intermediate field of $\KK/\bbk(Z)$, the last group $\Gal(\KK/\EE)$ is a subgroup
of $\Gal(\KK/\bbk(Z))=\Gal_\pi(\bbk)$.

Now suppose that
$\FF$ contains the algebraic closure of $\bbk$ and $Z$ is a rational variety, so that  $\bbk(Z)$ is a purely
transcendental extension of $\bbk$. 
Let $\defcolor{\LL}\vcentcolon=\bbk\cap\KK$ be the subfield of elements that are algebraic over
$\bbk$, which is a Galois extension of $\bbk$. 
Then in~\eqref{Eq:fieldExtensions}
$\EE=\LL(Z)$ is the compositum of $\LL$ and $\bbk(Z)$ and therefore  $\LL(Z)$ is algebraic and Galois
over $\bbk(Z)$.
Thus the Galois group $\Gal(\KK/\LL(Z))$ of $\KK$ over $\LL(Z)$ is a normal subgroup of $\Gal_\pi(\bbk)$.
Recalling that $\EE=\LL(Z)$ completes the proof.
\end{proof}

A consequence of Theorem~\ref{Th:GaloisUnderExtension} is that when $Z$ is rational, $\bbk=\QQ$,
and $\FF=\CC$ (or $\overline{\QQ}$), so that $\LL=\CC\cap\KK$, then  $\Gal_\pi(\CC)$ is a normal subgroup of
$\Gal_\pi(\QQ)$.
For the sparse polynomial systems of Section~\ref{S:GGsparse}, we often have that $\Gal_\pi(\CC)\neq\Gal_\pi(\QQ)$.
In all other cases that we know, these two groups are equal, but we lack a proof of that observation in general.

\subsection{Complex branched covers}\label{SS:complex}

Suppose that $\pi\colon X\to Z$ is a branched cover of complex varieties.
The
locus $U\subset Z$
where $\pi$ is proper and \'etale
is the open subset that is maximal with respect to inclusion such that the restriction
$\pi\colon\pi^{-1}(U)\to U$ is a covering space.
We will call $U$ the set of regular values of $\pi$.

The monodromy group $\Mon_\pi$ as defined in Section~\ref{SS:Galois_Monodromy} agrees with the usual notion of the monodromy group of the covering space
 \[
   \pi\ \colon\ \pi^{-1}(U)\ \longrightarrow\ U\,.
 \]
This is the group of permutations of a fiber $\pi^{-1}(z)$ obtained by lifting loops in $U$ that are based at $z$ to paths in
  $\pi^{-1}(U)$ that connect points in the fiber.
If $d$ is the degree of $\pi$, lifting based loops in $U$ to paths in a component $X'$ of $X_Z^{(d)}$ gives this equality.
For more on  covering spaces and monodromy groups, see~\cite{Hatcher,Munkres}.

The complement of any (Zariski) open subset $V$ of $Z$ has real codimension at least 2.
The loops in $U$ that generate the monodromy group can be chosen to lie in $V$ (by a change of base point if necessary).
A consequence is that the monodromy group $\Mon_\pi$ is equal to the monodromy group of any restriction
$\pi\colon\pi^{-1}(V)\to V$ to a Zariski open set $V$ such that this map is a covering space.

\subsection{Enriched Galois groups}

As Harris showed~\cite{Harris}, many  enumerative problems have Galois groups that are the full symmetric group $S_d$ on
their solutions.
We call such a Galois group/enumerative problem \demph{fully symmetric}.
It is a standard part of the Algebra curriculum that any finite group may arise as the Galois group of a branched cover.
Nevertheless, determining the possible Galois groups of a given class of enumerative problems (the \demph{inverse Galois problem} for
that class), as well as the Galois group of any particular enumerative problem is an interesting problem that is largely open.

Many techniques to study Galois groups in enumerative geometry are able to show that the Galois group $\Gal_\pi$ is either
$S_d$ or contains its subgroup \defcolor{$A_d$} of alternating permutations.
We call such an enumerative problem/Galois group  \demph{giant}.
While many enumerative Galois groups are 
known to be giant, we know of no natural enumerative problem whose
Galois group is the alternating group (besides those similar to $x\mapsto x^3$).

As we saw in the problem of 27 lines, when a Galois group fails to be fully symmetric, we expect there is a geometric reason for this
failure. 
That is, the set of solutions is enriched with extra structure that prevents the Galois group from being fully symmetric.
Consequently, we will call a Galois group or enumerative problem \demph{enriched} if its Galois group is not fully symmetric.

Let us recall some aspects of permutation groups.
A permutation group of degree $d$ is a subgroup $G$ of $S_d$.
Thus $G$ has a natural action on the set $\defcolor{[d]}\vcentcolon=\{1,\dotsc,d\}$, as well as on the subsets of
$[d]$. 
The group is \demph{transitive} if for any $i,j\in[d]$, there is an element $g\in G$ with $g(i)=j$.
More generally, for any $1\leq s\leq d$, $G$ is \demph{$s$-transitive} if for any distinct 
$i_1,\dotsc,i_s\in [d]$ and distinct $j_1,\dotsc,j_s\in[d]$, there is an element $g\in G$ with $g(i_m)=j_m$ for $m=1,\dotsc,s$.
That is, $G$ is $s$-transitive when it acts transitively on the set of distinct $s$-tuples of elements of $[d]$.
This has the following consequence.

\begin{proposition}\label{P:higher_transitivity}
  The monodromy group $\Mon_\pi$ of a branched cover $\pi\colon X\to Z$ is $s$-transitive if and only if the variety $X_Z^{(s)}$ is
  irreducible. 
\end{proposition}

Let $G$ be a transitive permutation group of degree $d$.
A \demph{block} of $G$ is a subset $B\subset [d]$ such that for every $g\in G$, either $gB=B$ or $gB\cap B=\emptyset$.
The subsets $\emptyset$, $[d]$, and every singleton are blocks of every permutation group.
If these trivial blocks are the only blocks, then $G$ is \demph{primitive} and otherwise it is \demph{imprimitive}. 

The Galois group $E_6$ for the problem of 27 lines is primitive, but it is not 2-transitive.
For the latter, observe that some pairs of lines on a cubic surface meet, while other pairs are disjoint.
These incidences provide an obstruction to 2-transitivity.

When $G$ is imprimitive, we have a factorization $d=ab$ with $1<a,b<d$ and there is a bijection
$[a]\times[b]\leftrightarrow[d]$ such that $G$ preserves the projection $[a]\times[b]\to [b]$.
That is, the fibers $\{[a]\times\{i\}\mid i\in [b]\}$ are blocks of $G$ and its action on this set of blocks gives a
homomorphism $G\to S_b$ with transitive image.
In particular, $G$ is a subgroup of the group of permutations of $[d]=[a]\times[b]$ which preserve the fibers of the projection
  $[a]\times[b]\to[b]$.
  This group is
the wreath product $S_a\wr S_b$, which is the semi-direct product $(S_a)^b\rtimes S_b$, where $S_b$ acts on $(S_a)^b$ by permuting
factors.

Imprimitivity has a geometric manifestation.
A branched cover $\pi\colon X\to Z$ is \demph{decomposable} if there is a nonempty Zariski open subset $V\subset Z$ and a
variety $Y$  such that $\pi$ factors over $V$,
 \begin{equation}\label{Eq:Ndecomposition}
  \pi^{-1}(V)\ \xrightarrow{\; \varphi\;}\ Y\ \xrightarrow{\; \psi\;}\ V\,,
 \end{equation}
with $\varphi$ and $\psi$ both nontrivial branched covers.
The fibers of $\varphi$ over points of $\psi^{-1}(v)$ are blocks of the action of $\Gal_\pi$ on $\pi^{-1}(v)$, which
implies that $\Gal_\pi$ is imprimitive. 
Pirola and Schlesinger~\cite{PirolaSchlesinger} observed that decomposability of $\pi$ is equivalent to imprimitivity of
$\Gal_\pi$.

\begin{proposition}\label{P:DecomposableIsImprimitive}
 A branched cover is decomposable if and only if its Galois group is imprimitive.
\end{proposition}  

Harris's geometric method~\cite{Harris} to show that a Galois group of an enumerative problem over $\CC$ is fully symmetric
involves two steps. 
First, show that $X_Z^{(2)}$ is irreducible, so that $\Mon_\pi$ is 2-transitive.
Next, identify an instance of the enumerative problem (a point $z\in Z$) with $d{-}1$ solutions, where exactly one
solution has multiplicity 2. 
This implies that a small loop in $Z$ around $z$ induces a simple transposition in $\Mon_\pi$.
This implies that $\Mon_\pi=S_d$, as $S_d$ is its only 2-transitive subgroup containing a simple transposition.
Jordan~\cite{J1870} gave a useful generalization of this last fact about $S_d$, which we use in Section~\ref{Sec:Computing}.

\begin{proposition}
 \label{P:Jordan}
 Suppose that $G\subset S_d$ is a permutation group.
   If $G$ is primitive and contains a $p$-cycle for some prime number
   $p<d{-}2$, then $G$  is giant.

 If  $G$ contains a $d$-cycle, a $d{-}1$-cycle, and a $p$-cycle for some prime number $p<d{-}2$ then $G=S_d$.
\end{proposition}

The first statement is the form of Jordan's Theorem found in, for example~\cite[Thm.\ 13.9]{Wielandt}.
It implies the second, as a permutation group containing both a $d$ and a $d{-}1$-cycle is primitive and also not contained
in the alternating group.

\section{Numerical Algebraic Geometry}  \label{Sec:NAG}
Methods from numerical analysis underlie algorithms that readily solve systems of polynomial equations
over $\CC$.
Numerical algebraic geometry uses this to represent and study algebraic varieties on a computer.
We sketch some of its fundamental algorithms, which will later be used for studying Galois groups.
For a more complete survey, see~\cite{NNLA}.

\subsection{Homotopy continuation}
When $\bbk=\CC$, solutions to enumerative problems, fibers of branched covers, and monodromy are all effectively computed
using algorithms based on numerical homotopy continuation.
This begins with a  \demph{homotopy}, which is a family $\calH(x;t)$ of systems of polynomials that interpolate between
the systems at $t=0$ and $t=1$ in a particular way:
We require that the variety $\calV(\calH(x;t))\subset \CC^n_x\times\CC_t$ contains a curve \defcolor{$C$} that is the
union of the 1-dimensional irreducible components of $\calV(\calH)$ which project dominantly to $\CC_t$.
We further require that $1\in\CC_t$ is a regular value of the projection $\pi\colon C\to\CC_t$, that $\pi$ is proper near 
$1$, and that $\calV(\calH(x;t))$ is smooth at all points of the fiber $\defcolor{W}\vcentcolon=\pi^{-1}(1)$.
The \demph{start system} is $\calH(x;1)$ and write \defcolor{$W$} for its set of isolated solutions, which we assume are
known. 
The \demph{target system} is $\calH(x;0)$ and we wish to use $\calH$ to compute the isolated solutions to the target
system.

Given a homotopy $\calH(x;t)$, we restrict $C$ to the points above an arc $\gamma\subset\CC_t$ with endpoints $\{0,1\}$
such that $\gamma$ avoids the critical values of $\pi\colon C\to\CC_t$ and points where $\pi$ is not proper, except
possibly at $t=0$. 
In what follows, we will take $\gamma$ to be the interval $[0,1]$, for simplicity.
This restriction is a collection of arcs in $C$, one for each point of $W$, which start at points of $W$ at $t=1$ and lie
above $(0,1]$. 
Some arcs may be unbounded near $t=0$, while the rest end in points of $\pi^{-1}(0)$, and all points of $\pi^{-1}(0)$ are
reached. 
Beginning with the (known) points of $W$, standard path-tracking algorithms~\cite{AG03} from numerical analysis may be used
to follow these arcs and 
compute the points of $\pi^{-1}(0)$.
When $\pi\colon C\to\CC_t$ is proper near $t=0$ and smooth above $t=0$, there are $|W|$ points in $\pi^{-1}(0)$ so that
each path gives a point of $\pi^{-1}(0)$.
In this case, the homotopy is \demph{optimal}.
Optimality is only one measure of complexity of solving systems.
  It is relevant as for many families of systems, the only known practical homotopy algorithms do not have that
   $\pi\colon C\to\CC_t$ is proper near $t=0$ and in fact follow exponentially many paths that diverge as $t\to 0$.
For more on numerical homotopy continuation, see~\cite{Mor87,SW05}.

The most straightforward optimal homotopy is a \demph{parameter homotopy}~\cite{LSY89,MS89}, in which the structure and
number of solutions 
of the start, target, and intermediate systems are the same.
A source for parameter homotopies is a branched cover $X\to Z$, where $Z$ is a rational variety and $X$ is a subvariety of
$\CC^n\times Z$.
Suppose that $f\colon \CC_t\to Z$ is a smooth rational curve with $f(0)$ and $f(1)$ lying in the open set $U$ of regular values
of $X\to Z$. 
Pulling back $X\to Z$ along $f$ gives a dominant map $\pi\colon \defcolor{f^*(X)}\to\CC_t$ with the same degree $d$ as
$X\to Z$. 
A generating set $\calH(x;t)$ of the ideal of $f^*(X)\subset \CC^n\times\CC_t$ gives a homotopy that is optimal as there
are $d$ solutions to both the start and target systems.

For example, suppose that $X\to Z=\PP^{19}$ is the branched cover~\eqref{Eq:Cubic_Family} from the problem of 27 lines.
Given smooth cubics $f_1$ and $f_0$, the pencil $\defcolor{f(t)}\vcentcolon=t f_1+(1-t)f_0$ is a map $\CC_t\to\PP^{19}$ as
above. 
A general line $\ell$ in $\PP^3$ is the span of points $[x_1,x_2,1,0]$ and $[x_3,x_4,0,1]$, for $(x_1,x_2,x_3,x_4)\in\CC^4$.
A general point on $\ell$ has the form $[ux_1+x_3,xu_2+x_4,u,1]$, for $u\in\CC$, and $\ell$ lies on the cubic $\calV(f(t))$
when $f(t)(ux_1+x_3,xu_2+x_4,u,1)$ is identically zero.
Thus, if we expand $f(t)(ux_1+x_3,xu_2+x_4,u,1)$ as a polynomial in $u$, the four coefficients of the resulting cubic are equations in
$x_1,\dotsc,x_4,t$ for the general line  $\ell$ to lie on the cubic $\calV(f(t))$.
Let $\calH(x;t)$ be these four coefficients.
When $\calV(f_1)$ has 27 lines of the given form, this is a homotopy, and if we knew the coordinates of those 27 lines, numerical homotopy
continuation using $\calH(x;t)$ could be used to compute the lines on $\calV(f_0)$.

\subsection{Witness sets}\label{SS:Witness}
Numerical homotopy continuation enables the reliable computation of solutions to systems of polynomial equations.
Numerical algebraic geometry uses this ability to solve as a basis for algorithms that study and manipulate varieties on a computer.
Its starting point is a witness set, which is a data structure for varieties in $\CC^n$~\cite{INAG,NAG,Witness}.
Suppose that $X\subset\CC^n$ is a union of irreducible components of the same dimension $m$ of a variety $\calV(F)$, where $F$ is a system
of polynomials.
If $L\subset\CC^n$ is a general linear subspace of codimension $m$, then $\defcolor{W}\vcentcolon=X\cap L$ is a transverse intersection consisting of
$\deg(X)$ points, called a \demph{linear section} of $X$.
The triple $(W,F,L)$ is a \demph{witness set} for $X$ (typically, $L$ is represented by $m$ linear forms).

Given a witness set $(W,F,L)$ for $X$ and a general codimension $m$ linear subspace $L'$, we may compute the linear section $W'=X\cap L'$
and obtain another witness set $(W',F,L')$ for $X$ as follows.
Let $L(t)\vcentcolon=t L + (1-t) L'$ be the convex combination of (the equations for) $L$ and $L'$, and form the homotopy $\calH(x;t)\vcentcolon=(F,L(t))$.
Path-tracking using $\calH(x;t)$ starting from the points of $W$ at $t=1$ will compute the points of $W'$ at $t=0$.
This instance of the parameter homotopy is called ``moving the witness set''.

Suppose that we have a third codimension $m$ linear subspace $L''$.
We may then use $W'$ to compute the linear section $W''=X\cap L''$, and then use $W''$ to return to $W$.
The arcs connect every point $w\in W$ to a point $w'\in W'$, then to a point $w''\in W''$, and finally to a possibly different point
$\sigma(w)\in W$.
This defines a permutation $\sigma$ of $W$.
The four points, as they are connected by smooth arcs, lie in the same irreducible component of $X$.
Thus the cycles in the permutation $\sigma$ refine the partition of $W$ given by the irreducible components of $X$.
Repeating this procedure with
possibly different linear subspaces $L',L''$, and then applying the trace test~\cite{traceLRS,SVW_trace},
leads to a \demph{numerical irreducible decomposition} of $X$; that is,
it computes the partition $W=W_1\sqcup\dotsb\sqcup W_r$, where $X_1,\dotsc,X_r$ are the irreducible components of $X$ and
$W_i\vcentcolon=X_i\cap L$.
This algorithm was developed in~\cite{SVW_decomposition,SVW_monodromy,SVW_trace}.

Several freely available software packages have implementations of the basic algorithms of Numerical Algebraic Geometry.
These include Macaulay 2~\cite{M2} in its Numerical Algebraic Geometry package~\cite{NAG4M2}, in Bertini~\cite{Bertini},
and in  HomotopyContinuation.jl~\cite{BT}.

\section{Fano Problems}\label{Sec:Fano}

Debarre and Manivel determined the dimension and degree of the variety of $r$-planes lying on general complete intersections in
$\PP^n$.
When this is zero-dimensional it is called a \demph{Fano problem}.
For example, the problem of 27 lines is a Fano problem.
Galois groups of Fano problems were studied classically by Jordan and Harris and recently by Hashimoto and Kadets, who nearly determined
the Galois group for each Fano problem.

\subsection{Combinatorics of Fano Problems}

Let \defcolor{$\GG(r,\PP^n)$} be the Grassmann variety defined over the complex numbers of $r$-dimensional linear subspaces of $\PP^n$,
which has dimension $(r{+}1)(n{-}r)$. 
Given a variety $X\subseteq\PP^n$, its \demph{Fano scheme} is the subscheme of $\GG(r,\PP^n)$ of $r$-planes lying on $X$.

Fano schemes may be studied uniformly when $X\subset\PP^n$ is a complete intersection.
For this, let $ \defcolor{d_\bullet}\vcentcolon= (d_1,\dotsc,d_s)$ be a weakly increasing list of integers greater than 1.
Suppose that  $F=(f_1,\dotsc,f_s)$ are homogeneous polynomials on $\PP^n$ with $f_i$ of degree $d_i$.
Let \defcolor{$\calV_r(F)$} be the Fano scheme of $r$-planes in $\calV(F)$.

Just as $\calV(F)$ has expected dimension $n{-}s$, there is an expected dimension for $\calV_r(F)$.
Let $f$ be a form on $\PP^n$ of degree $d$.
Its restriction to $H\in\GG(r,\PP^n)$ is a form of degree $d$ on $H$; as the dimension of the vector space of such forms is
$\binom{d+r}{r}$, we expect this to be the codimension of $\calV_r(f)$ in $\GG(r,\PP^n)$.
Thus the expected dimension of $\calV_r(F)$ is
\[
   \defcolor{\delta} =  \defcolor{\delta(r,n,d_\bullet)}\ \vcentcolon=\ (r{+}1)(n{-}r)\ -\ \sum_{i=1}^s\binom{d_i+r}{r}\,.
\]

Write  \defcolor{$\CC^{(r,n,d_\bullet)}$} for the vector space of homogeneous polynomials $F=(f_1,\dotsc,f_s)$ in $n{+}1$
  variables with $f_i$ of degree $d_i$.
Debarre and Manivel~\cite{Debarre} showed that there is a dense open subset $U=\defcolor{U_{(r,n,d_\bullet)}}\subset\CC^{(r,n,d_\bullet)}$
with the following property: 
For $F\in U$, if $\delta\ge 0$ and $n-s\ge 2r$, then  $\calV_r(F)$ is a smooth variety of dimension $\delta$, and if
$\delta<0$ or $n{-}s<2r$, then $\calV_r(F)$ is empty.
A \demph{Fano problem} is the enumerative problem of determining $\calV_r(F)$ for $F\in U_{(r,n,d_\bullet)}$, when
$\delta(r,n,d_\bullet)=0$ and $n{-}s\ge 2r$.

Since the Grassmannian has Picard group generated by $O(1)$ induced by its Pl\"ucker embedding, when  $\delta\ge 0$ and
$n{-}s\ge 2r$ and $F\in U$, the Fano variety  $\calV_r(F)$ has a well-defined degree.
Standard techniques in intersection theory allow this degree to be computed, using that $\calV_r(F)$ is the vanishing
of sections of appropriate vector bundles on $\GG(r,\PP^n)$.
(These are $\mbox{Sym}_{d_i}(T)$, where $T$ is the dual of the tautological $(r{+}1)$-subbundle on the Grassmannian.)

This leads to a formula for this degree.
For that,
define the polynomials
\[
    \defcolor{Q_{r,d}(x)}\ =\ \prod_{a_0 + \dotsb + a_r = d}(a_0 x_0 + \dotsb + a_r x_r)\ \in\  \ZZ[x_0,\dotsc,x_r]\ 
  \quad a_i\in\mathbb{Z}_{\ge 0}
\]
as well as $\defcolor{Q_{r,d_\bullet}} = Q_{r,d_1}(x)\dotsb Q_{r,d_s}(x)$ and the Vandermonde polynomial
\[
   \mbox{V}_r(x)\ =\ \prod_{0\le i<j\le r}(x_i - x_j)\,.
\]
When $\delta(r,n,d_\bullet)=0$, $n-s\ge 2r$, and $F\in U_{(r,n,d_\bullet)}$, the degree \defcolor{$\deg(r,n,d_\bullet)$} of $\calV_r(F)$ is
the coefficient of $x_0^n x_1^{n-1}\dotsb x_r^{n-r}$ in the product $Q_{r,d_\bullet}(x)\mbox{V}_r(x)$ \cite[Thm.\ 4.3]{Debarre}.
Table~\ref{SmallFano} gives these degrees for all Fano problems with a small number of solutions.
Here, $D_n$ refers to the
Coxeter group
of order $2^{n-1}n!$.
\newcommand{\SP}{\mbox{\hspace{4pt}}}
\begin{table}[htb]
  \caption{Small Fano problems.}
  \label{SmallFano}
  \def\arraystretch{1.2}
  \begin{tabular}{||c|c|c|c|c||}
    \hline
    \SP$r$\SP & \SP$n$\SP & $d_\bullet$ & $\#$ of solutions & Galois Group\\
    \hline\hline
    1 & 4 & $(2,2)$ & 16 & $D_5$\\
    \hline
    1 & 3 & $(3)$ & 27 & $E_6$\\
    \hline
    2 & 6 & $(2,2)$ & 64 & $D_7$\\
    \hline
    3 & 8 & $(2,2)$ & 256 & $D_9$\\
   \hline
    4 & 10 & $(2,2)$ & 1024 & $D_{11}$\\
    \hline
    1 & 4 & $(5)$ & 2875 & $S_{2875}$ \\
   \hline
   5 & 12 & $(2,2)$ & 2096 & $D_{13}$\\
   \hline
   $\vdots$ & $\vdots$ & $\vdots$ &$\vdots$&$\vdots$\\
   \hline
    \hline
  \end{tabular}
\end{table}

\subsection{Galois groups of Fano problems}
Consider the incidence correspondence,
 \[
   \raisebox{-30pt}{\begin{picture}(230,63)(-1,0)
    \put(0,52){$\Gamma\ \vcentcolon=\ \{(F,H)\in \CC^{(r,n,d_\bullet)}\times\GG(r,\PP^n)\,\mid\, F|_H = 0\}$\,.}
    \put(2,48){\vector(0,-1){31}}  \put(6,32){\small$\pi$}
    \put(-1,2){$\CC^{(r,n,d_\bullet)}$}
   \end{picture}}
 \]
The fiber over a general complete intersection $F\in U_{(r,n,d_\bullet)}$ is the Fano variety $\calV_r(F)$.
When we have a Fano problem,
$\pi$ is a branched cover of degree $\deg(r,n,d_\bullet)$.
We define the Galois group of the Fano problem to be $\defcolor{\Gal_{(r,n,d_\bullet)} = \Gal_\pi}$.  

The study of Galois groups of Fano problems began with Jordan \cite{J1870} with the problem of 27 lines on a smooth cubic surface, which
has data $(1,3,(3))$.
By observing the incidence structure of the lines on a smooth cubic, Jordan determined that the Galois group over $\CC$ is a subgroup of
$E_6$, $\Gal_{(1,3,(3))}\subseteq E_6$. 

Harris~\cite{Harris} showed that Jordan's inclusion is an equality, $\Gal_{(1,3,(3))} = E_6$, and then generalized this,
showing that $\Gal_{(1,n,(2n-3))}$ is fully symmetric for $n\ge 4$.
For this, he used the interpretation of the Galois group as a monodromy group.
Using arguments from algebraic geometry, when $n\ge 4$ he showed that the monodromy group is 2-transitive and contains a
simple transposition.

Hashimoto and Kadets~\cite{HK}
revisited this topic, determining these groups in many cases.
They first considered Fano problems of linear spaces on the
intersection of two quadrics.

\begin{proposition}\label{P:HK-1}
  All Fano problems of $r$-planes on the intersection of two quadrics in $\PP^{2r+2}$ are enriched, and
  \[
      \Gal_{(r,2r+2,(2,2))}\ =\ D_{2r+3}\,.
  \]
\end{proposition}

That the Fano problems of lines on space cubics and linear spaces on intersections of two quadrics
are enriched may be understood in that they are the only Fano problems where the $r$-planes on $\calV(F)$
are expected to intersect. 
As in the problem of 27 lines, the generic incidence structure prevents the Galois group from being fully symmetric.
In all other cases, they showed that the Galois group is giant.

\begin{proposition}\label{P:HK-2}
  Any Fano problem that is not of the form $(1,3,(3))$ or $(r,2r{+}2,(2,2))$, has giant Galois group.
\end{proposition}

(Recall that Harris showed that the Fano problems $(1,k,(2k{-}3))$ for $k>3$ are all fully symmetric.)
In Section~\ref{Sec:ComputingNumerically}, we describe a method based on numerical homotopy continuation which computes
monodromy permutations with high probability, when $\bbk=\CC$.
Using methods of numerical certification, Yahl~\cite{Yahl} extended Harris's methods to prove that several Fano problems of
moderate size have full symmetric Galois groups.
He proved the following.

\begin{theorem}\label{Th:Yahl}
  For each Fano problem $(r,n,d_\bullet)$ not equal to $(1,3,(3))$ or $(r,2r{+}2,(2,2))$ for $r\ge 1$ and with
  fewer than 75,000 solutions, the Fano Galois group $\Gal_{(r,n,d_\bullet)}$ is the full symmetric group.
\end{theorem}

This result is the product of determining Galois groups for the Fano problems of Table~\ref{BigFano},
each of which was shown to be the full symmetric group.
New ideas are needed to settle whether or not the remaining Fano problems are full symmetric.

\begin{table}[htb]
  \caption{Moderate-sized finite Fano problems.}
  \label{BigFano}
  \def\arraystretch{1.2}
  \begin{tabular}{||c|c|c|c||}
    \hline
    \SP$r$\SP  & $n$ & $d_\bullet$ & \SP$\deg(r,n,d_\bullet)$\SP\\
    \hline\hline
    1 & 7 & $(2,2,2,2)$ & 512\\
    \hline
    1 & 6 & $(2,2,3)$ & 720\\
    \hline
    2 & 8 & $(2,2,2)$ & 1024\\
    \hline
    1 & 5 & (3,3) & 1053\\
    \hline
    1 & 5 & (2,4) & 1280\\
    \hline
    1 &\SP 10\SP &\SP (2,2,2,2,2,2)\SP & 20480\\
    \hline
    1 & 9 & (2,2,2,2,3) & 27648\\
    \hline
    2 &\SP 10\SP & (2,2,2,2) & 32768\\
    \hline
    1 & 8 & (2,2,3,3) & 37584\\
    \hline
    1 & 8 & (2,2,2,4) & 47104\\
    \hline
    1 & 7 & (3,3,3) & 51759\\
    \hline
    1 & 7 & (2,3,4) & 64512\\
    \hline
  \end{tabular}
\end{table}

\section{Galois groups of sparse polynomial equations}
\label{S:GGsparse}

We work over the complex numbers. 
With modifications due to separability and constants (e.g. $\EE$ in proof of Theorem~\ref{Th:GaloisUnderExtension}),
much of this holds over an arbitrary field.

The Bernstein-Kushnirenko Theorem gives an upper bound on the number of solutions in the algebraic torus $(\CC^\times)^n$
to a system of polynomials.
This bound depends on the monomials which appear in the equations (their support).
When the equations are general given their support, this bound is attained.
The family of all systems with a given support forms a branched cover and therefore has a Galois group.
Esterov identified two structures in the support which imply that the Galois group is imprimitive, and showed that if they
are not present, then the Galois group is full symmetric.
It remains an open problem to determine the Galois group when it is imprimitive.

\subsection{Systems of sparse polynomial equations}

A \demph{(Laurent) monomial} in $n$ variables $x_1,\dotsc,x_n$ with exponent vector $\alpha = (\alpha_1,\dotsc,\alpha_n)\in\ZZ^n$ is
 \[
    \defcolor{x^\alpha}\ \vcentcolon=\ x_1^{\alpha_1}x_2^{\alpha_2}\dotsb x_n^{\alpha_n}\,.
 \]
This is a character of the algebraic torus $(\CC^\times)^n$. 
A \demph{(Laurent) polynomial} $f$ over $\CC$ is a linear combination of monomials,
\[
  f\ =\ \sum c_\alpha x^\alpha\qquad c_\alpha\in \CC\,.
\]
For a nonempty finite set $\calA\subseteq\ZZ^n$,
if the  above sum is restricted to $\alpha\in\calA$, then $f$ has \demph{support} $\calA$.
  Write \defcolor{$\CC^\calA$} for the set of polynomials $f$ with support $\calA$.

Given a collection $\defcolor{\calAdot}=(\calA_1,\calA_2,\dotsc,\calA_n)$ of nonempty finite subsets of
$\ZZ^n$, write 
 \[
 \defcolor{\CC^{\calAdot}}\  \vcentcolon=\  \CC^{\calA_1}\times\CC^{\calA_2}\times\dotsb\times \CC^{\calA_n}
 \]
for the vector space of $n$-tuples $F=(f_1,\dotsc,f_n)$ of polynomials, where $f_i$ has support $\calA_i$, for each $i=1,\dotsc,n$.
An element $F\in \CC^{\calAdot}$ is a square system of polynomials whose solutions are those $x\in(\CC^\times)^n$ such that 
 \[
   f_1(x_1,\dotsc,x_n)\ =\ 
   f_2(x_1,\dotsc,x_n)\ =\ \dotsb\ =\ 
   f_n(x_1,\dotsc,x_n)\ =\ 0\,,
 \]
written $F(x) = 0$.
We call $F$ a \demph{sparse polynomial system with support $\calAdot$}.

Given supports $\calAdot = (\calA_1,\dotsc,\calA_n)$, define the incidence variety
\[
   \raisebox{-30pt}{\begin{picture}(230,63)(-1,0)
       \put(0,52){$\Gamma \ =\
       \defcolor{\Gamma_{\calAdot}}\ \vcentcolon=\
  \left\{(F,x)\in \CC^{\calAdot}\times(\CC^\times)^n \mid F(x) = 0 \right\}\,$.}
    \put(2,48){\vector(0,-1){31}}  \put(6,32){\small$\pi$}
    \put(-4,2){$\CC^{\calAdot}$}
   \end{picture}}
 \]
It is equipped with projections $\pi\colon\Gamma\to  \CC^{\calAdot}$ and $p\colon\Gamma\to( \CC^\times)^n$.
The fiber $p^{-1}(x)$ for $x\in(\CC^\times)^n$ is the set of all polynomials $(f_1,\dotsc,f_n)$ with $f_i(x)=0$ for each $i$.
Observing that $f_i(x)=0$ is a non-zero linear equation on $\CC^{\calA_i}$, we see that
$p^{-1}(x)\subset \CC^{\calAdot}$ is the product of $n$ hyperplanes and thus has codimension $n$.
Consequently, $\Gamma\to (\CC^\times)^n$ is a vector bundle, and therefore irreducible, and it has dimension equal to $\dim \CC^{\calAdot}$.

Thus the map $\pi\colon\Gamma\to \CC^{\calAdot}$ is a branched cover when $\pi$ is dominant, equivalently, when a generic system
$F\in \CC^{\calAdot}$ has a positive, finite number of solutions in $(\CC^\times)^n$.
The number of solutions to a generic system is determined by the polyhedral geometry of its support, which we review. For convex bodies $K_1,\dotsc,K_n\subset\RR^n$ and nonnegative real numbers, $t_1,\dotsc,t_n\in\RR_{\geq0}$, Minkowski proved that the volume
of the Minkowski sum 
 \[
   t_1K_1+\dotsb+t_nK_n\ \vcentcolon=\ \{ t_1 x_1+\dotsb + t_nx_n \mid x_i\in K_i\}
 \]
is a homogeneous polynomial of degree $n$ in $t_1,\dotsc,t_n$.
Its coefficient of $t_1\dotsb t_n$ is the \demph{mixed volume} of $K_1,\dotsc,K_n$.
For supports $\calAdot = (\calA_1,\dotsc,\calA_n)$, let \defcolor{$\MV(\calAdot)$} be the mixed
volume of their convex hulls, $\conv(\calA_1),\dotsc,\conv(\calA_n)$.
This is described in detail in \cite{Ewald}.
We state the Bernstein-Kushnirenko Theorem~\cite{Bernstein75,Kushnirenko2}.

\begin{theorem}[Bernstein-Kushnirenko]
  A system $F\in \CC^{\calAdot}$ has at most $\MV(\calAdot)$ isolated solutions in $(\CC^\times)^n$.
  This bound is sharp and is attained for generic $F\in \CC^{\calAdot}$.
\end{theorem}

Thus $\pi\colon \Gamma_{\calAdot}\to\CC^{\calAdot}$ is a branched cover of degree $\MV(\calAdot)$ if and only if
$\MV(\calAdot)> 0$, which Minkowski determined as follows.
For a nonempty subset $I\subseteq\defcolor{[n]}\vcentcolon=\{1,\dotsc,n\}$, write $\defcolor{\calA_I}\vcentcolon=(\calA_i\mid i\in I)$ and
let \defcolor{$\ZZ\calA_I$} be the affine span of the supports in $\calA_I$.
This is the free abelian group generated by all differences $\alpha{-}\beta$ for $\alpha,\beta\in \calA_i$ for
some $i\in I$. 
Then $\MV(\calAdot) = 0$ if and only if there exists a 
subset $I\subseteq[n]$ such that $|I|>\rank(\ZZ\calA_I)$.
One direction is obvious.
When  $|I|>\rank(\ZZ\calA_I)=m$, then there is a change of variables so that the
subsystem of polynomials with indices in $I$ has more equations than variables.
In particular, $\MV(\calAdot)\neq 0$ implies that $\defcolor{\ZZ\calAdot}\vcentcolon=\ZZ\calA_{[n]}$ has full rank $n$.

\subsection{Galois groups of sparse polynomial systems}

Suppose that $\calAdot$ is a collection of supports with $\MV(A_\bullet)>0$.
Write \defcolor{$\Gal_{\calAdot}$} for the Galois group of the corresponding branched cover
$\pi\colon\Gamma_{\calAdot}\to \CC^{\calAdot}$.
Esterov~\cite{Esterov} studied these groups, identifying two structures which imply that $\Gal_{\calAdot}$ is imprimitive.

\begin{example}\label{Ex:Lacunary}
  Let $n=1$ and suppose that we have a univariate polynomial $f(x)$ of the form $g(x^3)$, for $g$ a univariate polynomial
  with $g(0)\neq 0$.
  Observe that $\ZZ\calAdot\subset 3\ZZ$.
  The zeroes of $f$ ($\{x\in\CC\mid f(x)=0\}$) are cube roots of the zeroes of $g$, and  the group of cubic roots of unity
  acts freely on the zeroes of $f$.
  These orbits are blocks of the action of the Galois group of $f$.
  When $g$ has two or more roots, there is more than one orbit, and the action of the Galois group is imprimitive.

  For another example, suppose that our system $F$ is
  \[
  1-2xy+3x^2-5xy^{-1}+7x^2y^2\ =\
  1+2y^2-4xy+8xy^3+16y^4\ =\ 0\,.
  \]
  This may be written as $G(xy,x/y)$, where $G(s,t)$ is
  \[
  1-2s +3st - 5t + 7s^2\ =\
  1+2s/t -4 s + 8 s^2/t + 16 s^2/t^2\ =\ 0\,.
  \]
  Given any of the (five) solutions $(s^*,t^*)$ of $G=0$, there are two solutions to $F$,
  \begin{equation}\label{Eq:Blocks}
    [\pm\bigl(\sqrt{s^*t^*}\,,\ \sqrt{s^*/t^*}\;\bigr)\,.
  \end{equation}
  (The branches of the two quare roots are chosen coherently, so that $s^*=\sqrt{s^*t^*} \sqrt{s^*/t^*}$ and
  $t^*=\sqrt{s^*t^*}\sqrt{s^*/t^*}$.)
  These pairs~\eqref{Eq:Blocks} are blocks of the action of the Galois group, showing that it is imprimitive.\hfill$\diamond$
\end{example}

Generalizing this, call $\calAdot$ \demph{lacunary} if $\ZZ\calAdot\neq\ZZ^n$.
If $\MV(\calAdot)>[\ZZ^n\colon \ZZ\calAdot]$, then it is \demph{strictly lacunary}.

\begin{example}\label{Ex:Triangular}
  Suppose now that our system $F$ is 
  \begin{eqnarray*}
    f(y)&=&1+2y+5y^2+ 2y^3 +y^4\\
    g(y,x)&=& 1+ 2y-3z+ 5y^2 -8 yz +13 z^2 + 21yz^2-34z^3\,.
  \end{eqnarray*}
  We first find the four roots $y^*$ of $f(y)=0$.
  For each root $y^*$, the polynomial $g(y^*,z)$ is a cubic with three solutions.
  Thus the Galois group of $F$ permutes the roots of $f$ with the roots of $g$ for each forming blocks for its
  action.
  Observe that $F$ has the form $f(y)=g(y,z)=0$.\hfill$\diamond$
\end{example}

This example generalizes.
Call $\calAdot$ \demph{triangular} if there exists a nonempty proper subset $I\subsetneq [n]$ such that
$\rank(\ZZ\calA_I) = |I|$. 
In this case there is a change of variables so that the equations indexed by $I$ involve only the first $|I|$ variables.
We write $\MV(\calA_I)$ for the mixed volume of the supports of the polynomials indexed by $I$ as polynomials in the first
$|I|$ variables.
We say $\calAdot$ is \demph{strictly triangular} if $1<\MV(\calA_I)<\MV(\calAdot)$.

We provide more details in Section~\ref{Sec:Applications}.
This is explained fully in~\cite{Esterov}, and algorithmically in~\cite[Sect.\ 2.3]{DSS}.
When $\calAdot$ is neither lacunary nor strictly triangular, Esterov showed that $\Gal_{\calAdot}$ is 2-transitive.
Then, he showed that a small loop around the discriminant of these systems generates a simple transposition,
which shows that $\Gal_{\calAdot}$ is full symmetric. 

\begin{theorem}[Esterov]
  \label{Th:Esterov}
  Let $\calAdot$ be a set of supports with $\MV(\calAdot)>0$.
  If $\calAdot$ is neither lacunary nor strictly triangular, then $\Gal_{\calAdot}$ is the full symmetric group.
  If $\calAdot$ is strictly lacunary or strictly triangular, then $\Gal_{\calAdot}$ is imprimitive.
  If $\calAdot$ is lacunary but not strictly lacunary, then  $\Gal_{\calAdot}$ is the group $\Hom(\ZZ^n/\ZZ\calAdot,\CC^\times)$ of roots of
  unity. 
\end{theorem}

When $\calAdot$ is either strictly lacunary or strictly triangular, Esterov's theorem does not determine the group $\Gal_{\calAdot}$ explicitly.
As it is imprimitive, the Galois group $\Gal_{\calAdot}$ is a subgroup of a certain wreath product.
It may be a proper subgroup, as the following example shows.

\begin{example}\label{Ex:subgroupOfWreathProduct}
  Let $n=2$ and suppose that $\calA$ is consists of the vertices of the $2\times 2$ square and its center point $(1,1)$, which we show
  below. 
 \[
  \begin{picture}(102,77)(-18,-1)
     \put(0,0){\includegraphics{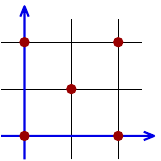}}
     \put(-18,64){\small$(2,0)$}     \put(62,64){\small$(2,2)$}
     \put(-18,-1){\small$(2,0)$}     \put(62,-1){\small$(2,0)$}
  \end{picture}
 \]
 Let $\calAdot\vcentcolon=(\calA,\calA)$.
 Its mixed volume is $\MV(\calAdot)=8$, which is twice the area of the square.
 Thus a general system of polynomials with support $\calAdot$ has eight solutions in $(\CC^\times)^n$.
 The lattice $\ZZ\calAdot$ has index 2 in $\ZZ^2$, so solutions come in four pairs of symmetric points, $(x,y)$ and
 $(-x,-y)$. 
 These pairs are preserved by the Galois group, showing that it is a subgroup of the wreath product
 $S_2\wr S_4$.
 It can be shown that $\Gal_{\calAdot} = (S_2\wr S_4)\cap A_8$ and is thus a proper subgroup of this wreath
 product.\hfill$\diamond$ 
\end{example}

This example is due to Esterov and Lang~\cite{EsterovLang}, who gave conditions which imply that the Galois group is the full wreath product,
for certain lacunary systems.
Despite this, there is no known criteria for when that occurs, not even a conjecture about which groups
may occur as Galois groups of sparse polynomial systems~\cite{EsterovLang_new,Yahl_GG}.
Also, it is not clear what can be said about Galois groups of sparse polynomials over other fields than the complex
numbers.

\section{Computing Galois Groups}
\label{Sec:Computing}

Understanding Galois groups of enumerative problems has both benefited from and inspired the
development of and use of computational tools.
We discuss an adaptation of the well-known symbolic method of computing cycle types of Frobenius elements and then several
methods based on numerical homotopy continuation.
For a prime $p\in\ZZ$, write $\defcolor{\FF_p}\vcentcolon=\ZZ/p\ZZ$ for the field with $p$ elements.

\subsection{Galois groups of univariate polynomials}
  \label{S:univariate}
 Determining Galois groups of univariate polynomials is an old and challenging problem in mathematics.
 Modern algorithms, such as Stauduhar's~\cite{Stauduhar} and its more recent improvements~\cite{FiKl}
 are implemented in systems incuding MAGMA~\cite{MAGMA}.
 These are effective for rational polynomials $f(t)$ of moderate degree $d$.

 These may be used to study Galois groups of branched covers $\pi\colon X\to Z$ defined over $\QQ$ when $Z$ is rational.
 We may formulate a fiber $\pi^{-1}(z)$ for $z\in Z$ as a system $F(x;z)$ of rational polynomials in some number, $n$, of
 variables and parameters $z\in Z$.
 For $z\in Z(\QQ)$ this is a system of rational polynomials in $x$.
 Methods based on Gr\"obner bases may be used to compute a rational univariate representative~\cite{RUR}, which is a
 univariate polynomial $f\in \QQ[t]$ of degree $d=\deg_\pi$ whose roots generate the field $\KK$ of definition of the fiber
 $\pi^{-1}(z)$. 
 
 Once computed, $f(t)$ becomes an input to these algorithms to compute $\Gal(\KK/\QQ)\subset\Gal_\pi(\QQ)$.
 Experience suggests that this inclusion is often an equality.
 A drawback to this approach is that these polynomials typically have large coefficient size (arithmetic height).

 For an example, beginning with a cubic surface in $\PP^3$ whose coefficients are random two-digit decimal integers, the
 resulting system $F$ has only four variables $x$.
 It takes 30 seconds to compute $f(t)$, which has degree $27$ and coefficients that have on average 165 digits, and then MAGMA
 takes seconds to tell us that its Galois group is $E_6$.
 The study of reality in~\cite{RSSS} involved computing these representatives $f(t)$.
 It was restricted to $d\lesssim 30$ and $n\lesssim 10$ varibles as it was infeasible to compute $f(t)$ for larger values
 of $d$ and $n$.

 The obstruction is the arithmetic height of the polynomials encountered.
 A recent project is studying enriched problems on the Lagrangian Grassmannian~\cite{LGn}.
 For one problem of degree $d=16$ in $n=5$ variables, a representative was computed whose decimal integer coefficients 
 had over 3500 digits.
 For larger enriched problems in this preliminary study, with $d\in\{16,24,32,64,128\}$  and $n\sim 9$ variables,
 computing $f(t)\in\QQ[t]$ is infeasible.
 Alternative methods for studying Galois groups in enumerative geometry that avoid this problem of arithmetic height
 are developed in subsequent sections. 
 
\subsection{Frobenius elements}

Let $f\in\ZZ[x]$ be a monic irreducible univariate polynomial
with splitting field $\KK$, a finite Galois extension of $\QQ$.
Let $\defcolor{\calO}\subset \KK$ be the ring of integers in $\KK$, the set of elements of $\KK$ that are integral over
$\ZZ$. 
For every prime $p\in\ZZ$ not dividing the discriminant of $f$, there is a \demph{Frobenius element} 
$\defcolor{\sigma_p}\in\Gal(\KK/\QQ)$
(well-defined up to conjugacy) 
in the Galois group of $\KK$ over $\QQ$ that restricts to the Frobenius 
automorphism above $p$:
For every prime ideal \defcolor{$\varpi$} of $\calO$ with $\varpi \cap \ZZ = \langle p \rangle$
($\varpi$ is \demph{above} $p$),
and every $z\in\calO$, we have $\sigma_p(z)\equiv z^p \mod\varpi$.
That is, $\sigma_p$ restricts to the Frobenius automorphism on $\calO/\varpi$.
If $f$ is not monic, then we first invert the primes dividing the leading coefficient of $f$.
The existence of such Frobenius elements is explained in~\cite[\S\S~VII.2]{Lang}.

The cycle type of a Frobenius element $\sigma_p$ (as a permutation of the roots of $f$) is given by the degrees
of the irreducible factors in $\FF_p[x]$ of $\defcolor{f_p}\vcentcolon=f\mod p$, as these factors give prime ideals $\varpi$ above $p$.
Indeed, if $g$ is an irreducible factor of $f$ of degree $r$ with corresponding prime ideal $\varpi$,
then $\calO/\varpi\simeq \FF_p[x]/\langle g\rangle$ is a finite field with $p^r$ elements.
The Frobenius automorphism on $\FF_p[x]/\langle g\rangle$ acts on the roots of $g$ as a cycle of length $r$.
The cycle type of $\sigma_p$  records how $p$ splits in
$\calO$ and is also called the \demph{splitting type} of $\sigma_p$ or of $f_p$. 
The prime $p$ does not divide the discriminant exactly when $f_p$ is squarefree and it has the same degree as $f$.
This gives an algorithm to compute cycle types of Frobenius elements of $\Gal(\KK/\QQ)$:
For a prime $p$ with $\deg(f_p)=\deg(f)$, factor the reduction $f_p$, and if no factor is repeated, record the degrees of
the factors.

This method is particularly effective due to the Chebotarev Density Theorem~\cite{StLe,Cheb}:
Let $\KK/\QQ$ be a  Galois extension and $\lambda$ a cycle type of an element in $\Gal(\KK/\QQ)$.
Define  \defcolor{$n_\lambda$} to be the fraction of elements in $\Gal(\KK/\QQ)$ with cycle type $\lambda$.
Then the density of primes $p\leq N$ such that the Frobenius element $\sigma_p$ has splitting type $\lambda$ tends to $n_\lambda$ as
$N\to\infty$.
Loosely speaking, for $p$ sufficiently large, Frobenius elements are distributed uniformly in $\Gal(\KK/\QQ)$.

Table~\ref{T:S_6} illustrates this for 
%
%
$f=x^6-503x^5-544x^4-69x^3-152x^2-49x-763$,
\begin{table}[htb]
\caption{Frobenius elements for $f$.}\label{T:S_6}
 \centering
$\begin{array}{||r|r|r|r|r|r|r|r|r|r|r||}\hline
     1^6 & 1^4,2 & 1^2,2^2 & 2^3 & 1^3,3 & 1, 2, 3 & 3^2 & 1^2, 4 & 2, 4 & 1, 5 &  6 \rule{0pt}{12pt} \\\hline\hline
     1   &   15  &   45   &  15 &   40  &  120   &  40 &   90   &  90  &  144 & 120 \\\hline\hline
     3   &   12  &   24   &   9 &   47  &  146   &  32 &  112   &  71  &  121 & 143 \\\hline
   .989 & 15.02 & 44.97 & 14.99 & 40.07 & 120.03 & 39.95 & 89.87 & 89.97 & 144.24 & 119.9 \\\hline\hline
\end{array}$
\end{table}
which has Galois group $S_6$.
%
%
The headers in the first row are the cycle types (conjugacy classes) of permutations in $S_6$, expressed using the frequency
representation for cycle type in which $(2^3)$ indicates three 2-cycles.
The second row contains the sizes of each conjugacy class.
The third row records how many of the first $720=6!$ primes $p$ not dividing the
discriminant\footnote{$18972006774677773002386748159696=2^4\cdot 3^{12}\cdot 7\cdot 29\cdot 2633\cdot 88805021 \cdot 47006055979$.}
did $f_p$ have the corresponding splitting type. 
For the last row, we repeated this calculation for the first 
$720\cdot 10^5$ primes larger than $10^8$.
We display the observed number that had a given splitting type, divided by $10^5$ for comparison.
The convergence guaranteed by the Chebotarev Density Theorem is evident.

Determining the splitting type of Frobenius elements gives information about Galois groups, including information about the
distribution of cycle types in a Galois group when sufficiently many are computed.
For example, if the Galois group $\Gal$ is known to be a subgroup of a particular permutation group $G$, knowing the cycle
types of relatively few elements often suffices to show that $\Gal=G$, as Proposition~\ref{P:Jordan} does when $G$ is the
symmetric group. 
If we do not have a candidate for $\Gal$, then computing many Frobenius elements may help to predict the Galois group with
a high degree of confidence, by the Chebotarev Density Theorem.
Both approaches were crucial for the computations of Schubert Galois groups in Section~\ref{S:SchubertGaloisGroups}.

\subsection{Frobenius elements for branched covers}
Frobenius elements are also a tool for studying Galois groups in enumerative geometry using symbolic computation.
Suppose that $\pi\colon X\to Z$ is a branched cover of degree $d$ of irreducible affine varieties defined over $\ZZ$, and
that $Z$ is a smooth, rational variety.
Restricting to an open subset of $Z$, we may also assume that the map $\pi\colon X\to Z$ is proper and \'etale.
By the results in~\cite[\S\S~VII.2]{Lang}, for each prime $p\in\ZZ$ and closed point $z_p\in Z(\FF_p)$, there is a Frobenius
element $\defcolor{\sigma_{z_p}}\in\Gal_\pi(\QQ)$ as before.
Its cycle type is the splitting type of the fiber $\pi^{-1}(z_p)$.

Given a prime $p$ and a cycle type $\lambda$ of an element in the Galois group $\Gal_\pi(\QQ)$, we may consider the density
of points $z\in Z(\FF_p)$  such that the corresponding Frobenius element has conjugacy class $\lambda$.
Ekedahl~\cite{Ekedahl} showed that in the limit as $p\to\infty$, this density tends to $n_\lambda$, the density of the
conjugacy class in $\Gal_\pi(\QQ)$.

This theoretical result gives an algorithm to study $\Gal_\pi(\QQ)$ (we touched on this in~\S~\ref{S:univariate}).
To begin, we recast the derivation of a Frobenius element in more elementary terms.
If $z\in Z(\ZZ)$ is an integer point, it is in particular a rational point and the fiber  $\pi^{-1}(z)$ is a reduced
zero-dimensional subscheme of $X$.
As $X$ is affine,  we may assume that $X\subset\mathbb{A}^r$.
If we extend scalars to the algebraic closure $\overline{\QQ}$ of $\QQ$,  $\pi^{-1}(z)$ consists of $d$
points in $\overline{\QQ}\,^r$.
Their coordinates generate a subfield \defcolor{$\KK_z$} that is a Galois extension of $\QQ$ whose Galois group is a
subgroup of $\Gal_\pi(\QQ)$.
For all except finitely many prime numbers $p$, both $z$ and $\pi^{-1}(z)$ have reductions \defcolor{$z_p$} and
\defcolor{$\pi^{-1}(z_p)$} modulo $p$, and thus a Frobenius element $\defcolor{\sigma_{z_p}}\in\Gal(\KK_z/\QQ)$.
This is conjugate to the Frobenius element of the first paragraph.

Assume that $\pi\colon X\to Z$ is a branched cover  of irreducible varieties defined over $\ZZ$ with $Z$ an open
subset of an affine space $\bAA^m(\ZZ)$.
All enumerative problems we discuss have this form, as $Z$ is typically a variety of parameters (coefficients of
polynomials or entries of matrices representing flags).
Replacing $X$ by an open subset, we have that $X\subset \bAA^m(\ZZ)\times\bAA^n(\ZZ)$ is an affine variety with ideal
$I\subset\ZZ[z,x]$. 
Specializing $I$ at an integer point $z\in Z(\ZZ)$ gives the ideal \defcolor{$I(z)$} of the fiber $\pi^{-1}(z)$.
The splitting type of the fiber at $p$ may be determined by a primary decomposition of the ideal $I(z)$ modulo $p$.
(This may also be accomplished with a rational univariate representative.)
This allows us to sample from $\Gal_\pi(\QQ)$ and is a tool for studying this Galois group.
We illustrate this method for lines on cubic surfaces.

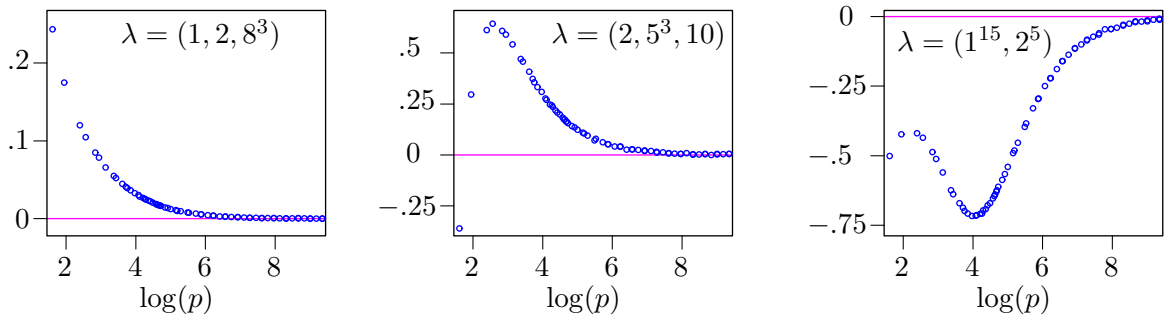
\begin{figure}[htb]
 \centering
  \begin{picture}(135,113)(-10,-21)
    \put(0,0){\includegraphics{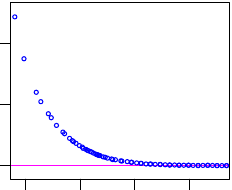}}
    \put(-10,68){\small$.2$}
    \put(-10,38){\small$.1$}
    \put(-7, 8){\small$0$}
    \put(33,78){\small$\lambda=(1,2,8^3)$}
    \put(9,-10){\small$2$} \put(35.7,-10){\small$4$} \put(62.01,-10){\small$6$} \put(88.2,-10){\small$8$}
    \put(40,-22){\small$\log_e(p)$}
  \end{picture}
  \begin{picture}(151,113)(-25,-21)
    \put(0,0){\includegraphics{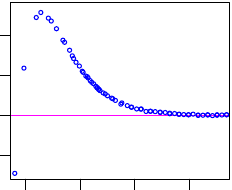}}
    \put(-16,73){\small$.5$}
    \put(-16,53){\small$.25$}
    \put(-13,33){\small$0$}
    \put(-24,13){\small$-.25$}
    \put(42,78){\small$\lambda=(2,5^3,10)$}
    \put(9,-10){\small$2$} \put(35.7,-10){\small$4$} \put(62.01,-10){\small$6$} \put(88.2,-10){\small$8$}
    \put(40,-22){\small$\log_e(p)$}
  \end{picture}\quad
  \begin{picture}(145,113)(-25,-21)
    \put(0,0){\includegraphics{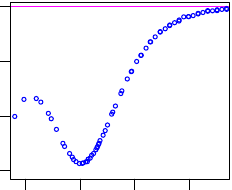}}
    \put(-12,85){\small$0$}
    \put(-24,58.67){\small$-.25$}
    \put(-24,32.33){\small$-.5$}
    \put(-24, 6){\small$-.75$}
    \put(10,75){\small$\lambda=(1^{15},2^5)$}
    \put(9,-10){\small$2$} \put(35.7,-10){\small$4$} \put(62.01,-10){\small$6$} \put(88.2,-10){\small$8$}
    \put(40,-22){\small$\log_e(p)$}
   \end{picture}
   
   \caption{Relative discrepancy, $\frac{E_{p,\lambda}}{n_\lambda}-1$, against $\log_{e}(p)$ at 69 primes $p$,
   for splitting types $\lambda\in\{(1,2,8^3), (2,5^3,10), (1^{15},2^5)\}$. }
   \label{F:27_plots} 
\end{figure}

Consider the branched cover $\pi\colon\Gamma\to \PP^{19}$ of lines on cubic surfaces~\eqref{Eq:Cubic_Family}.
For each of 69 primes $p$ between  $5$ and $11579$, we used Singular~\cite{Singular} to determine the splitting type of the
27 lines for many (70 to 220 million) randomly 
chosen smooth cubic surfaces in $\PP^{19}(\FF_p)$, and compared that to the distribution of cycle types in the Galois group $E_6$.
Recall that $n_\lambda$ is the density of elements in $E_6$ with cycle type $\lambda$.
For a prime $p$, let \defcolor{$E_{p,\lambda}$} be the empirical density, the observed fraction of surfaces whose lines had splitting type
$\lambda$.
By Ekedahl's Theorem, $\lim_{p\to\infty}E_{p,\lambda}=n_\lambda$.
(This same limit holds if we replace smooth cubics in $\PP^{19}(\FF_p)$ by isomorphism classes of smooth cubics~\cite{BFL}.)
Figure~\ref{F:27_plots} presents some data from our calculation.
We observe that the convergence in Ekedahl's Theorem may be slow.
The full computation is archived on the web
page\footnote{{\tt https://FrankSottile.github.io/research/stories/27\_Frobenius/}}.

There are algorithms to compute this decomposition implemented in software such as Macaulay2~\cite{M2} or
Singular~\cite{Singular}.
While these rely on Gr\"obner bases, they can be
unreasonably effective, as they also take advantage of fast Gr\"obner basis calculation in 
positive characteristic, for example the F5 algorithm~\cite{F5}.

\subsection{Computing Galois groups numerically}\label{Sec:ComputingNumerically} 

In Section~\ref{SS:Witness}, we discussed how moving a witness set $W$ to another one, $W'$, to a third, $W''$, and then
back to $W$ computes a permutation $\sigma$ of $W$.
This is readily adapted to computing a permutation of a fiber of a branched cover, which is an element of its Galois group
$\Gal_\pi$. 
While computing several such monodromy permutations only gives a subgroup of $\Gal_\pi$, that may be sufficient to
determine it~\cite{LS09}. 
We explain a numerical method from~\cite{ngalois} that computes a generating set for the Galois group and another
numerical method to study transitivity.

Given a branched cover $\pi\colon X\to Z$ of degree $d$ over $\CC$ with $Z$ rational, let $U\subset Z$ be the regular
locus so that $\pi^{-1}(U)\to U$ is a covering space.
Suppose that we have computed all points in a fiber $\pi^{-1}(z)$ for some point $z\in U$.
Choosing other points $z',z''\in U$, we may construct three parameter homotopies that move the points of $\pi^{-1}(z)$ to those of
$\pi^{-1}(z')$, to $\pi^{-1}(z'')$, and then back to $\pi^{-1}(z)$.
The tracked paths give a permutation $\sigma$ of the fiber $\pi^{-1}(z)$.

Writing the points of $\pi^{-1}(z)$ in some order $(w_1,\dotsc,w_d)$ gives a point in the fiber of $X^{(d)}_Z$ over $z$.
(Recall that $X^{(d)}_Z$ was used in Section~\ref{SS:Galois_Monodromy} to give a geometric construction of $\Gal_\pi$.)
The $d$-tuples of computed paths between the fibers $\pi^{-1}(z)$, $\pi^{-1}(z')$, $\pi^{-1}(z'')$, and back to
$\pi^{-1}(z)$ give a path in $X^{(d)}_Z$ from the point $(w_1,\dotsc,w_d)$ to the point $(\sigma(w_1),\dotsc,\sigma(w_d))$
in the same fiber. 
These paths all lie in the same irreducible component of $X^{(d)}_Z$, showing that $\sigma\in\Gal_\pi$.

This may be used to compute many permutations in $\Gal_\pi$, giving an increasing sequence of subgroups of $\Gal_\pi$.
As with computing Frobenius elements, this may suffice to determine $\Gal_\pi$.
For example, if the computed subgroup of $\Gal_\pi$ is $S_d$, then $\Gal_\pi=S_d$ is full-symmetric.
This method was used in~\cite{LS09} to show that several Schubert Galois groups (see Section~\ref{S:SchubertGaloisGroups})
were full-symmetric, including one with $d=17589$.
In that computation, every time a new permutation $\pi$ was found, GAP~\cite{GAP4} was called to test if the computed set of permutations
generated the symmetric group.

A drawback is that numerical path tracking may be inexact, which can lead to false conclusions (this is known as path-crossing).
A consequence of the calculation in~\cite{LS09} was the implementation of an algorithm~\cite{HS12} for {\it a posteriori} certification of
the computed solutions to a system of polynomials, based on Smale's $\alpha$-theory~\cite{S86}.
Certification using Krawczyk's method from interval arithmetic~\cite{Krawczyk} has also been implemented~\cite{BRT}, providing another
approach.
More substantially, algorithms were developed~\cite{BL12,BL13,DuLe} to certify path-tracking and thereby certifiably compute monodromy.

This approach of computing monodromy may be improved to compute a generating set of the Galois group~\cite{ngalois}.
Given a branched cover $\pi\colon X\to Z$ as above, restricting to an open subset of $Z$ and compactifying, we may assume that $Z=\PP^N$.
The branch locus of $\pi\colon X\to\PP^N$ is a hypersurface $B\subset\PP^N$.
Let $z\in \defcolor{U}\vcentcolon=\PP^n\smallsetminus B$.
Lifting loops in $U$ based at $z$ gives a surjective homomorphism from the fundamental group of $U$ to the monodromy group
of the cover $\pi^{-1}(U)\to U$~\cite{Hatcher,Munkres}.

A witness set for $B$ can be used to obtain a generating set for the fundamental group  of $U$.
Suppose that $\ell\cap B$ is a linear section of the hypersurface $B$, so that $\ell\simeq\PP^1$ is a line.
Zariski~\cite{Zar} showed that the inclusion $\ell\cap U \hookrightarrow U$ induces a surjection from the fundamental group of $\ell\cap U$
to the fundamental group of $U$.
As the fundamental group of $\ell\cap U$ is generated by based loops around each of the (finitely many) points of
$\ell\smallsetminus(B\cap \ell)$, lifts of these loops generate the Galois group $\Gal_\pi$ of the branched cover.

In~\cite{ngalois}, this method is demonstrated on the branched cover $\Gamma\to\PP^{19}$~\eqref{Eq:Cubic_Family} from the
problem of 27 lines.
The branch locus $B$ is the set of singular cubics, which forms a hypersurface on $\PP^{19}$ of degree 32, so that
a general line $\ell$ in $\PP^{19}$ meets $B$ transversally in 32 points.
The computed permutations for a particular choice of $\ell$ (given in Figure~5 in~\cite{ngalois})  generate $E_6$.
Each permutation is a product of six disjoint 2-cycles in $S_{27}$.
Here is one,
\[
   (1,6)(4,13)(8,25)(10,19)(11,16)(20,27)\ .
\]
That a loop around a point of $\ell\cap B$ gives a permutation that is the product of six disjoint 2-cycles is a manifestation of the
  enriched structure of this enumerative problem; Above a general point of $B$, there are six solutions of multiplicity 2.
  Contrast this with the result of Esterov~\cite{Esterov} from Section~\ref{S:GGsparse} where a single loop around $B$ gave a simple
  transposition.

Similar ideas were used to establish Theorem~\ref{Th:Yahl}, except that rather than compute a full witness set for the branch locus, a
  single point $z\in B$ of the branch locus was explicitly computed having the property that the fiber over $z$ consists of a single solution of multiplicity 2 and is otherwise smooth.
  Lifting a small loop around such $z$ gives a simple transposition---this was sufficient to show that those Fano problems were full symmetric.

We mention another method from~\cite{ngalois} involving transitivity.
Let $\pi\colon X\to Z$ be a branched cover of degree $d$.
By Proposition~\ref{P:higher_transitivity}, for any $1\leq s\leq d$, $s$-transitivity of the Galois group $\Gal_\pi$ is equivalent to  the
irreducibility of the variety $X_Z^{(s)}$.
Numerical irreducible decomposition may be used to determine the (ir)reducibility of $X_Z^{(s)}$ and thereby determine
whether or not $\Gal_\pi$  is $s$-transitive.
Details and an example involving the problem of 27 lines are given in~\cite[Sect.\ 4]{ngalois}.


\section{Galois groups in Schubert calculus} \label{S:SchubertGaloisGroups}

In his seminal book, ``Kalkul der abz{\"a}hlenden Geometrie''~\cite{Sch1879} Schubert presented methods for computing the
number of solutions to problems in enumerative geometry.
Justifying these methods was  Hilbert's 15th problem~\cite{HilbertE}, and they collectively came to be known as ``Schubert's
Calculus''. 
A central role was played by the Grassmannian and its Schubert cycles/varieties.
Schubert and others studied these objects further, and now Schubert varieties and the interplay of their geometry,
combinatorics, and algebra make them central objects in combinatorial algebraic geometry~\cite{Fu97} and other areas of
mathematics. 
This study is also called \demph{Schubert calculus}.
We are concerned with the overlap of these versions of Schubert calculus---problems in enumerative geometry that involve
intersections of Schubert varieties in Grassmannians and flag manifolds.

These \demph{Schubert problems} form a rich and well-understood class of examples that has long served as a laboratory for
investigating new phenomena in enumerative geometry~\cite{KL72}.
Thousands to millions of Schubert problems are computable and therefore may be studied on a computer.
Recently, this has also included reality in enumerative geometry~\cite{IHP}, and the resulting experimentation 
generated conjectures~\cite{HHMST,RSSS,So00} and examples~\cite{ModFour} concerning reality in Schubert calculus.
These in turn have helped to inspire proofs of some conjectures~\cite{EG02,EGSV,KP,LP,MT16,MTV_Sh,MTV_R,Purbhoo}.

Vakil's geometric Littlewood-Richardson rule~\cite{Va06a} gave a new tool~\cite{Va06b} for investigating Galois groups of
Schubert problems (\demph{Schubert Galois groups}) on Grassmannians.
He used it to discover an infinite family of Schubert problems on Grassmannians with enriched Galois groups.
The study of reality in flag manifolds uncovered another infinite family of enriched Schubert problems in manifolds of
partial flags~\cite[Thm.~2.18]{RSSS}. 
Subsequent results and constructions have led to the expectation that a Schubert Galois group in type $A$ should
be an iterated wreath product of symmetric groups, together with an understanding of the structure of enriched Schubert
problems.
It also appears that a Schubert Galois group should  not depend upon the base field.
Despite this, we are far from a classification, and the study has been largely 
limited to Grassmannians and type $A$ flag manifolds.
\begin{remark}\label{R:LGn}
  A preliminary study of enriched problems on the Lagrangian Grassmannian~\cite{LGn} has found a greater variety of
  Schubert Galois groups, including $(\ZZ/2\ZZ)^n$ and Schubert problems geometrically equivalent to the Fano 
  problems $(r,2r+2,(2,2))$.\hfill$\diamond$
\end{remark}

After describing Schubert problems in Grassmannians, in Section~\ref{SS:Enriched_SP} we construct Schubert problems whose
Galois groups (over any field) are symmetric groups $S_b$ acting on flags of subsets of $[b]\vcentcolon=\{1,\dotsc,b\}$;
This gives many enriched Schubert problems on flag manifolds in type $A$.
We also present a conjectural solution to the inverse Galois problem for Schubert calculus in type $A$.
In Section~\ref{SS:IGP_SP} we describe a general construction of Schubert problems whose Galois groups are expected to be
wreath products of two Schubert Galois groups. 
Our last section discusses results on Schubert Galois groups that are leading to an emerging picture of a
possible classification of Schubert problems in type $A$ by their Galois groups.

\subsection{Schubert problems}\label{SS:Scubert_Problems}

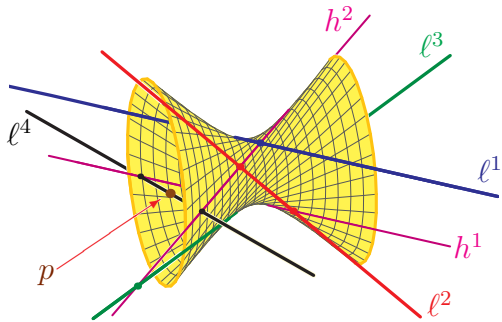
\begin{figure}[htb]
%
%
\centerline{
  \begin{picture}(189,108)(-1,5)
   \put(3,0){\includegraphics[height=4.1cm]{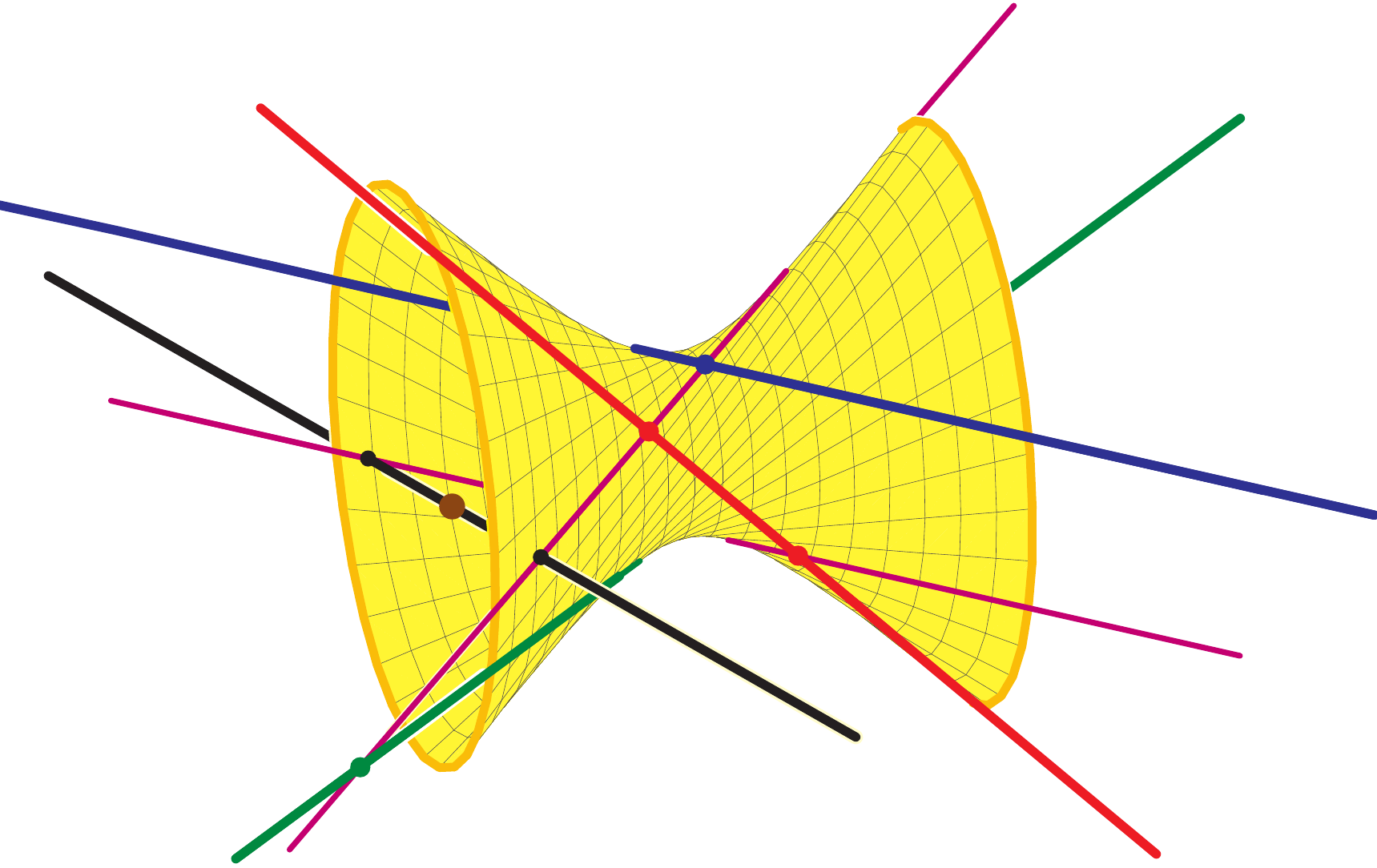}}
   \put(161,  2){{\color{Red}$\ell^2$}}
   \put(180, 53){{\color{Blue}$\ell^1$}}
   \put(158,100){{\color{Green}$\ell^3$}}
   \put(  2, 67){$\ell^4$}
   \put(171, 24){{\color{Magenta}$h^1$}}
   \put(122,110){{\color{Magenta}$h^2$}}
   \thicklines
   \put(21,20.1){{\color{White}\line(3,2){39}}}
   \thinlines
   \put(21,20.1){{\color{Red}\vector(3,2){39}}}
   \put(14, 16){{\color{Brown}$p$}}
  \end{picture}
}
 \caption{Problem of four lines.}\label{F:4lines}
\end{figure}

Consider the classical Schubert problem:  ``Which lines in $\PP^3$ meet each of four general lines?''
Three mutually skew lines   
{\color{Blue}$\ell^1$}, {\color{Red}$\ell^2$}, and {\color{Green}$\ell^3$} 
lie on a unique hyperboloid (Fig.~\ref{F:4lines}).
This hyperboloid has two rulings.
One contains {\color{Blue}$\ell^1$}, {\color{Red}$\ell^2$}, and {\color{Green}$\ell^3$}, 
and the second consists of the lines meeting them.
A general fourth line, $\ell^4$, meets the hyperboloid in two points,
and through each of these points there is a unique line in the second ruling.
These two lines, {\color{Magenta}$h^1$} and {\color{Magenta}$h^2$}, are the solutions to this instance of the problem of four lines.
Its Galois group is the symmetric group $S_2$:
Indeed, the solutions move as $\ell_4$ moves and rotating $\ell_4$ $180^\circ$ about the point $p$ will interchange the two lines.

More generally, a Schubert problem involves determining the linear subspaces of a vector space that have specified positions with
respect to certain fixed, but general linear subspaces.
For the problem of four lines, if we replace projective space $\PP^3$ by $\bbk^4$, the lines become 2-dimensional
linear subspaces.
Thus the problem of four lines is to determine the 2-planes in $\bbk^4$ that meet four given 2-planes nontrivially.

We introduce some terminology.
First, fix integers $1\leq m<n$.
The collection of all $m$-dimensional subspaces of $\bbk^n$ is the \demph{Grassmannian} $\defcolor{\Gr(m,n)}$ (also written
$\defcolor{\Gr(m,\bbk^n)}$), which is an algebraic manifold of dimension $m(n{-}m)$.
In Section~\ref{Sec:Fano}, this space was written $\GG(m{-}1,\PP^{n-1})$.

The Grassmannian has distinguished Schubert varieties.
These depend upon the choice of a (complete) \defcolor{flag}, which is a collection
$\defcolor{F}\colon F_1\subset F_2\subset\dotsb\subset F_n=\bbk^n$ of linear subspaces with $\dim{F_i}=i$.
Given a flag $F$, a Schubert variety is the collection of all $m$-planes having a given position with respect to $F$.
A position is encoded by a \demph{partition}, which is a weakly decreasing sequence of nonnegative integers
$\defcolor{\lambda}\colon \lambda_1\geq\dotsb\geq \lambda_m\geq 0$ with $\lambda_1\leq n{-}m$.
For a flag $F$ and partition $\lambda$, the corresponding Schubert variety is
 \begin{equation}\label{Eq:Schubert_variety}
   \defcolor{\Omega_\lambda F}\ \vcentcolon=\ \{H\in\Gr(m,n)\mid \dim H\cap F_{n-m+i-\lambda_i}\geq i\ \mbox{ for }i=1,\dotsc,m\}.
 \end{equation}
Setting $\defcolor{|\lambda|}\vcentcolon=\lambda_1+\dotsb+\lambda_m$, the Schubert variety $\Omega_\lambda F$ has codimension
$|\lambda|$ in $\Gr(m,n)$.

Only $m$ of the $n$ subspaces $F_i$ in $F$ appear in the definition~\eqref{Eq:Schubert_variety} of the Schubert variety $\Omega_\lambda F$. 
If $i<m$ and $\lambda_i = \lambda_{i+1}$, then the condition on $H$ in~\eqref{Eq:Schubert_variety} from
$\lambda_{i+1}$ implies the condition on $H$ for $\lambda_i$.
Those $i$ with $\lambda_i>\lambda_{i+1}$ (or $i=m$ with $\lambda_m>0$) are \demph{essential}.
When $(m,n)=(2,4)$ and $\lambda=(1,0)$, the essential condition is when $i=1$.
Indeed, $\Omega_{(1,0)}F=\{H\in\Gr(2,4)\mid \dim H\cap F_2\geq 1\}$.
In $\PP^3$, this is the set of lines $\PP H$ that meet the fixed line $\PP F_2$.

A \demph{Schubert problem} is a list $\defcolor{\lamDot}=\lambda^1,\dotsc,\lambda^s$ of partitions with $\sum_j|\lambda^j|=m(n-m)$, the
dimension of the Grassmannian.
An \demph{instance} of 
$\lamDot$ is given by a choice $\defcolor{\Fdot}=(F^1,\dotsc,F^s)$ of flags.
The solutions to this instance form the set of $m$-planes $H$ that have position $\lambda^j$ with respect to the flag $F^j$, for each $j$.
This set is the intersection 
 \begin{equation}\label{Eq:Schubert_Instance}
    \Omega_{\lambda^1}F^1\,\bigcap\,
    \Omega_{\lambda^2}F^2\,\bigcap\,\dotsb\,\bigcap\,
    \Omega_{\lambda^s}F^s\,.
 \end{equation}
Kleiman~\cite{KL74} showed that this intersection is transverse when the flags are general and $\bbk$ has
  characteristic zero.
  Transversality for positive characteristic with $\bbk=\overline{\bbk}$ is due to Vakil~\cite{Va06b}.
When the flags are general,
this implies that  for each solution (point $H$ in~\eqref{Eq:Schubert_Instance}), the inequalities in~\eqref{Eq:Schubert_variety}
for each pair $\lambda^j,F^j$ hold with equality.
Also, the number of solutions does not depend upon the (general) flags.
Write \defcolor{$d(\lamDot)$} for this number, which may be computed using algorithms in Schubert's calculus.

Let  \defcolor{$\Fl(n)$} be the space of complete flags in $\bbk^n$ and consider the incidence variety:
 \begin{equation}\label{Eq:Schubert_branched_cover}
   \raisebox{-30pt}{\begin{picture}(300,64)(-15,-1)
     \put(-3,52){$\defcolor{\Gamma_{\lamDot}} \vcentcolon= \{(H,F^1,\dotsc,F^s)\in\Gr(m,n)\times\Fl(n)^s\,\mid\,$}
                       \put(143,35){$H\in\Omega_{\lambda^i}F^i\ \mbox{for}\ i=1,\dotsc,s\}$}
    \put(2,47){\vector(0,-1){30}}  \put(6,31){\small$\pi_{\lamDot}$}
    \put(-15,2){$\Fl(n)^s$}
   \end{picture}}
 \end{equation}
The total space $\Gamma_{\lamDot}$ of this Schubert problem is irreducible, as it is a fiber bundle over the Grassmannian $\Gr(m,n)$ with
irreducible fibers (this is explained in~\cite[Sect.\ 2.2]{SW15}).
The fiber of $\pi_{\lamDot}$ over $(F^1,\dotsc,F^s)\in\Fl(n)^s$ is the intersection~\eqref{Eq:Schubert_Instance}.
Since this is transverse and consists of $d(\lamDot)$ points for general flags, $\pi_{\lamDot}$ is a branched cover of degree $d(\lamDot)$.
We write \defcolor{$\Gal_{\lamDot}$} for its Galois group, which we call a \demph{Schubert Galois group}.
We will often omit the field $\bbk$ as the constructions we give are independent of the field.

\subsection{Some enriched Schubert problems}\label{SS:Enriched_SP}
We present a construction of many enriched Schubert problems on Grassmannians and flag manifolds.
  These are based on the following generalization of
the problem of four lines:
Let $1<b$ be an integer and consider the 2-planes $h\subset\bbk^{2b}$ that meet each of four general
$b$-planes $K^1,\dotsc,K^4$ in at least a one-dimensional subspace.
Since $2b-2+1-(b-1)=b$, this Schubert problem is given by four partitions, each equal to $(b{-}1,0)$.

We explain how to solve this Schubert problem.
As the $K^i$ are general, we have $K^i\oplus K^j=\bbk^{2b}$ for $i\neq j$.
Then $K^3,K^4$ are graphs of isomorphisms $f_3,f_4\colon K^1\to K^2$, and
$\defcolor{\varphi}\vcentcolon=f_4^{-1}\circ f_3$ is a linear isomorphism of $K^1$, and any isomorphism may occur in this
way.  
If $\ell=\varphi(\ell)\subset K^1$ is a $\varphi$-stable line (one-dimensional subspace) in $K^1$, then
$f_3(\ell)=f_4(\ell)$ and $H\vcentcolon=\ell\oplus f_3(\ell)$ is a solution to this enumerative problem.
Furthermore, all solutions have this form, as $H\cap K^i$ for $i=1,\dotsc,4$ are four lines in the same 2-plane $H$. 
As the subspaces $K^i$ are general, the linear transformation $\varphi$ is semi-simple and therefore has
$b=\dim(K^1)$ distinct $\varphi$-stable lines.
Thus this Schubert problem has $b$ solutions (we are working over $\overline{\bbk}$ for these solutions).

Note that the monodromy group for the enumerative problem of stable lines of a semi-simple linear transformation $\varphi$ is
the full symmetric group $S_b$ acting on the set of 1-dimensional $\varphi$-stable linear subspaces of $K^1$.
In the notation of Section~\ref{SS:Galois_Monodromy}, $Z=GL(K^1)$ and
$X=\{(\varphi,\ell)\in Z\times\PP(K^1)\mid\varphi(\ell)\subset\ell\}$, and then
\[
  X^{(b)}_Z\ =\
  \{(\varphi,\ell_1,\dotsc,\ell_b)\in Z\times\PP(K^1)^b\mid\varphi(\ell_i)\subset\ell_i\ i\neq j\Rightarrow \ell_i\neq\ell_j\}\,,
\]
which is stable under the action of $S_b$ given by permuting the factors of $\PP(K^1)^b$.
The same construction shows that the Galois group of this Schubert problem is the full symmetric group $S_b$ acting 
naturally as permutations on the set $\defcolor{[b]}\vcentcolon=\{1,\dotsc,b\}$.
Using other means, Vakil~\cite[\S~3.14]{Va06b} also shows that the Galois group is $S_b$.

\begin{example}\label{Eq:Derksen-Vakil}
 Vakil~\cite[\S~3.14]{Va06b} used these problems in $\Gr(2,2b)$ to construct an infinite family of Schubert problems with enriched Galois
 groups.
 Let $1\leq a < b$ and consider the Schubert problem in $\Gr(2a,2b)$ of $2a$-planes that meet each of four general $b$-planes
 $K^1,\dotsc,K^4$ in at least an $a$-dimensional linear subspace.
 The previous argument generalizes:
 Let $f_3,f_4\colon K^1\to K^2$ and $\varphi\vcentcolon=f_4^{-1}\circ f_3\colon K^1\to K^1$ be the linear isomorphisms determined by $K^1,\dotsc,K^4$.
 Every solution has the form $L\oplus f_3(L)$ for $L=\varphi(L)\subset K^1$ a $\varphi$-stable $a$-dimensional linear subspace.
 Consequently, $L$ is spanned by $a$ linearly independent eigenvectors of $\varphi$.
 Thus this Schubert problem has $\binom{b}{a}$ solutions.

 The symmetric group $S_b$ acts naturally on the set \defcolor{$\binom{[b]}{a}$} of subsets of $[b]$ of cardinality $a$, and this
 argument shows that this permutation group (written \defcolor{$S\binom{[b]}{a}$}) is the Galois group of this Schubert problem.
 This action is not 2-transitive when $1<a< b{-}1$.
 It preserves the dimension of the intersection of a pair of solutions, and
 thus has at least $\min\{a,b{-}a\}$ distinct orbits on pairs of solutions. \hfill$\diamond$
\end{example}

We generalize Vakil's examples, while also generalizing~\cite[Thm.~2.18]{RSSS}.

\begin{example}\label{Eq:Enriched_flags}
  Suppose that $1\leq a_1<\dotsb<a_r<b$ are integers and write $\defcolor{\adot}$ for the sequence $a_1<\dotsb<a_r$.
  Let $\defcolor{\Fl(2\adot,2b)}$ be the space of partial flags of the form
 \[  
    F\ \colon\ F_{2a_1}\ \subset\ F_{2a_2}\ \subset\ \dotsb\ \subset\ F_{2a_r}\ \subset\ \bbk^{2b}\,,
 \]  
where $\dim F_{2a_i} = 2a_i$.
Fix four general $b$-planes $K^1,K^2,K^3,K^4$ in $\bbk^{2b}$.
Consider the Schubert problem that seeks the partial flags $F\in \Fl(2\adot,2b)$ such that for $i=1,\dotsc,r$,
\begin{equation}\label{Eq:flag_problem}
  \dim F_{2a_i}\bigcap K^j\ \geq\ a_i\qquad\mbox{ for all }\ j=1,\dotsc,4\,.
\end{equation}  
As before, $K^1,\dotsc,K^4$ give isomorphisms $f_3,f_4\colon K^1\to K^2$ and $\varphi=f_4^{-1}\circ f_3\colon K^1\to K^1$.
For each $1\leq i\leq r$, the solutions to~\eqref{Eq:flag_problem} are given by $L_{a_i}\oplus f_3(L_{a_i})$ where
$L_{a_i}\subset K^1$ is a $\varphi$-stable linear subspace of dimension $a_i$.

Consequently, the solutions to~\eqref{Eq:flag_problem} for all $i$ are in bijection with $\varphi$-stable flags
\[
     L_{a_1}\ \subset\ L_{a_2}\ \subset\ \dotsb\ \subset\ L_{a_r}\ \subset\ K^1\,,
\]
where $\dim L_{a_i}=a_i$.
Since $L_{a_i}$ is necessarily spanned by $a_i$ independent eigenvectors of $\varphi$, these are in bijection with flags of subsets of $[b]$:
\[
  \defcolor{\binom{[b]}{\adot}}\ \vcentcolon=\
  \{ T_1\subset T_2\subset\dotsb\subset T_r\subset[b]\,\mid\, |T_i|=a_i\}\,.
\]
Thus $\binom{[b]}{\adot}$ counts solutions to this Schubert problem and 
its Galois group is the symmetric group $S_b$, with its natural action on the set
$\binom{[b]}{\adot}$ of flags of subsets.
Write \defcolor{$S\binom{[b]}{\adot}$} for this permutation group.\hfill$\diamond$
\end{example}

This completes the following existence proof concerning Schubert Galois groups.

\begin{theorem}\label{Th:permutationSGG}
  For any positive integers $1\leq a_1<\dotsb<a_r<b$, there is a Schubert problem on the flag manifold $\Fl(2\adot,2b)$
  with Galois group $S\binom{[b]}{\adot}$.
\end{theorem}

These Schubert Galois groups form the basis for a conjectured solution to the Inverse Galois Problem in Schubert calculus.

\begin{conjecture}\label{C:IGP}
  A Galois group for a Schubert problem on a type A flag manifold is an iterated wreath product of permutation groups $S\binom{[b]}{\adot}$,
  and all such wreath products occur.

  Schubert Galois groups for Grassmannians are iterated wreath products of permutation groups $S\binom{[b]}{a}$,
  and all such wreath products occur.
\end{conjecture}

Conjecture~\ref{C:IGP} describes all known Schubert Galois groups---we discuss that and more in Section~\ref{S:SC_history}.
Additionally, all Schubert problems we know of with enriched Galois groups are either among those described in
Examples~\ref{Eq:Derksen-Vakil} or~\ref{Eq:Enriched_flags} or they are fibrations of Schubert problems, a structure we discuss in
Section~\ref{SS:IGP_SP} which is conjectured to give such wreath products.
Also, in all cases the Schubert Galois group does not depend upon the field $\bbk$.

\subsection{Compositions of Schubert problems}\label{SS:IGP_SP}
By Proposition~\ref{P:DecomposableIsImprimitive}, when a branched cover is decomposable, its Galois group is a subgroup of the wreath product
of the Galois groups of its factors. 
We explain how to construct a Schubert problem on a Grassmannian $\Gr(2a{+}m,2b{+}n)$ with decomposable branched cover.
This is built from one of the Schubert problems of Example~\ref{Eq:Derksen-Vakil} on $\Gr(2a,2b)$ and any Schubert problem $\muDot$ on
$\Gr(m,n)$ with $d(\muDot)>1$.
Conjecturally, its Galois group is the wreath product $(\Gal_{\muDot})^{\binom{b}{a}}\rtimes S\binom{[b]}{a}$.
This conjecture would establish existence in the Inverse Galois Problem for Schubert problems on Grassmannians.

It is convenient to represent a partition $\mu$ by its (Young) diagram, which is a left-justified array of boxes with $\mu_i$
boxes in row $i$.
Thus
\[
         (1)\ \longleftrightarrow\ \raisebox{-2pt}{\includegraphics{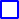}}\,, \qquad
       (2,2)\ \longleftrightarrow\ \raisebox{-6pt}{\includegraphics{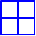}}\,, \quad\mbox{and}\quad\ 
   (3,2,1,1)\ \longleftrightarrow\ \raisebox{-14pt}{\includegraphics{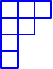}}\ .
\]
We omit any trailing 0s in a partition $\mu$.
Observe that the essential conditions in $\mu$ correspond to the boxes that form south east corners.
Consequently, a rectangular partition imposes a single incidence condition.

As the number $d(\muDot)$ of solutions to a Schubert problem $\muDot$ may be computed in the cohomology ring of the corresponding
Grassmannian~\cite{Fu97}, we often write a Schubert problem multiplicatively.
Thus $(\raisebox{-.5pt}{\includegraphics{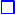}},\raisebox{-.5pt}{\includegraphics{pictures/1s}},
\raisebox{-.5pt}{\includegraphics{pictures/1s}},\raisebox{-.5pt}{\includegraphics{pictures/1s}})$,
which is the problem of four lines,  is also written
$\raisebox{-.5pt}{\includegraphics{pictures/1s}}\cdot\raisebox{-.5pt}{\includegraphics{pictures/1s}}\cdot
\raisebox{-.5pt}{\includegraphics{pictures/1s}}\cdot\raisebox{-.5pt}{\includegraphics{pictures/1s}}$ or as
$\raisebox{-.5pt}{\includegraphics{pictures/1s}}^4$.
The construction of a composition of Schubert problems is a bit technical, we will illustrate it first on the simplest example, when
  $\lamDot=\muDot$ are both the problem of four lines.

\begin{example}\label{Ex:First_Composition}
  Consider the Schubert problem $\kapDot$ in $\Gr(4,8)$ given by the partitions
 \begin{equation}\label{Eq:firstComposition}
   \kapDot\ =\ 
   \raisebox{-2pt}{\includegraphics{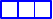}\,}\cdot\ 
   \raisebox{-2pt}{\includegraphics{pictures/3}\,}\cdot\ 
   \raisebox{-10pt}{\includegraphics{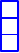}\,}\cdot\ 
   \raisebox{-10pt}{\includegraphics{pictures/111}\,}\cdot\ 
   \raisebox{-2pt}{\includegraphics{pictures/1}\,}\cdot\ 
   \raisebox{-2pt}{\includegraphics{pictures/1}\,}\cdot\ 
   \raisebox{-2pt}{\includegraphics{pictures/1}\,}\cdot\ 
   \raisebox{-2pt}{\includegraphics{pictures/1}\,}.
 \end{equation}
An instance of this Schubert problem is given by two 2-planes $\ell^1_2,\ell^2_2$, two 6-planes $L^3_6,L^4_6$, and four 4-planes
$K^1_4,\dotsc,K^4_4$, and its solutions are
 \begin{equation}\label{Eq:First_Composition}
   \{H\in\Gr(4,8)\,\mid\,
      \begin{array}{r} \dim H\cap \ell^i_2\geq 1\ \mbox{ for }i=1,2\\
                       \dim H\cap L^j_6\geq 3\ \mbox{ for }j=3,4\end{array} \ 
       \dim H\cap K^t_4\geq 1\ \mbox{ for }t=1,\dotsc,4
      \}\ .
 \end{equation}
Assume that these linear subspaces $\ell^i_2,L^j_6,K^t_4$ are in general position, which implies the dimension assertions that follow.
Let $\defcolor{\Lambda}\vcentcolon=\ell^1_2\oplus\ell^2_2\simeq \bbk^4$ and $\defcolor{M}\vcentcolon=L^3_6\cap L^4_6\simeq\bbk^4$.

If \defcolor{$H$} is a solution to~\eqref{Eq:First_Composition}, then $\dim H\cap\Lambda\geq 2$ and $\dim H\cap M\geq 2$.
As $\Lambda\cap M=\{0\}$ and $\dim H=4$, these inequalities are equalities.
Set $\defcolor{h}\vcentcolon=H\cap\Lambda\in G(2,\Lambda)$.
For $j=3,4$, the intersection $H\cap L^j_6$ has codimension 1 in $H$ and therefore $\dim h\cap L^j_6=1$.
Setting $\defcolor{\ell^j_2}\vcentcolon=\Lambda\cap L^j_6$, we see that $h$ is a solution to the instance of the problem of four lines given by
$\ell^1_2,\dotsc,\ell^4_2$.

Similarly, for each $i=1,\dotsc,4$, $H\cap M$ meets the 2-plane $\defcolor{k^i_2(h)}\vcentcolon=(h\oplus K^i_4)\cap M$.
Thus $H\cap M$ is a solution to the problem of four lines given by $k^1_2(h),\dotsc,k^4_2(h)$.
Conversely, given a solution $h\subset\Lambda$ to the problem of four lines given by $\ell^1_2,\dotsc,\ell^4_2$ and a solution
$h'\subset M$ to the  problem of four lines given by $k^1_2(h),\dotsc,k^4_2(h)$, their sum $h\oplus h'$ is a solution to the Schubert
problem~\eqref{Eq:First_Composition}. \hfill$\diamond$
\end{example}

Thus the branched cover $\pi\colon\Gamma_{\kapDot}\rightarrow \Fl(8)^8$ of this Schubert
problem~\eqref{Eq:Schubert_branched_cover} is decomposable. 
Indeed, let $U\subset\Fl(8)^8$ be the subset of flags in the general position used in
Example~\ref{Ex:First_Composition}. 
If we let $X\vcentcolon=\pi^{-1}(U)$ be the restriction of $\Gamma_{\kapDot}$ to this set of general instances, then
Example~\ref{Ex:First_Composition} shows that
we have a factorization $X\rightarrow Y \rightarrow U$~\eqref{Eq:decomposableFlags}.
Here, the fiber of $Y\to U$ over an instance in $U$ is the instance of $\raisebox{-.5pt}{\includegraphics{pictures/1s}}^4$
in $\Gr(2,\Lambda)$ given by $\ell_2^1,\dotsc,\ell_2^4$, and given a solution $h$ to this instance, the fiber of $X\to Y$
over $h$ is the instance of $\raisebox{-.5pt}{\includegraphics{pictures/1s}}^4$ in $\Gr(2,M)$ given by
$k^1_2(h),\dotsc,k^4_2(h)$.  

We make a definition inspired by this structure.

\begin{definition}\label{D:fiberedSchubertProblem}
  A Schubert problem $\kapDot$ is \demph{fibered} over a Schubert problem $\lamDot$ with fiber $\muDot$ if the branched
  cover $\Gamma_{\kapDot}\to\Fl(n)^s$ is decomposable, and it admits a decomposition
  \begin{equation}\label{Eq:decomposableFlags}
   X\ \longrightarrow\ Y\ \longrightarrow\ U \quad (U\ \subset\ \Fl(n)^s\ \mbox{ is open and dense})
  \end{equation}
  such that
  \begin{enumerate}
   \item
     fibers of $Y\to U$ are instances of $\lamDot$, 

   \item
     fibers of $X\to Y$ are instances of $\muDot$, and 

   \item
     general instances of $\lamDot$ and $\muDot$ occur as fibers in this way.
     
  \end{enumerate}
\end{definition}

We will call $\kapDot$ a \demph{fibration}.
This notion is developed in~\cite{GIVIX} and~\cite{SWY}, where the following is proven, which is a special case of~\cite[Lemma 15]{GIVIX}.

\begin{proposition}\label{P:wreath}
  If a Schubert problem $\kapDot$ is fibered over $\lamDot$ with fiber $\muDot$, then 
  $d(\kapDot)=d(\lamDot)\cdot d(\muDot)$, and its Galois group is a subgroup of the wreath product
  \[
     \Gal_{\kapDot}\ \subset\ \left( \Gal_{\muDot} \right)^{d(\lamDot)} \rtimes \Gal_{\lamDot}\,.
  \]  
\end{proposition}  

Consequently, the Schubert Galois group from Example~\ref{Ex:First_Composition} is a subgroup of the wreath product
$(S_2)^2\rtimes S_2$. 
In fact, its Galois group equals this wreath product~\cite[Sect.\ 5.5.2]{MSJ}.
This was shown by computing sufficiently many Frobenius elements.

The construction of  Example~\ref{Ex:First_Composition} was generalized in~\cite{SWY}.
Given two Schubert problems $\lamDot$ and $\muDot$ on possibly different Grassmannians, that paper describes how to use
them to build a new  Schubert problem \defcolor{$\lamDot\circ\muDot$} on another Grassmannian, called their
\demph{composition}. 
It uses combinatorics to prove that $d(\lamDot\circ\muDot)=d(\lamDot)\cdot d(\muDot)$,
and it is expected---but not proven---that  $\lamDot\circ\muDot$  is fibered over $\lamDot$ with fiber $\muDot$.

Next, it identifies a family of Schubert problems and shows that for any Schubert problem $\lamDot$ in that family, any composition
$\lamDot\circ\muDot$  is fibered over $\lamDot$ with fiber $\muDot$.
This  family includes all the Schubert problems of Example~\ref{Eq:Derksen-Vakil}.
We explain this construction when $\lambda$ is a Schubert problem of Example~\ref{Eq:Derksen-Vakil}, which is a motivation for
Conjecture~\ref{C:IGP}. 

Write \defcolor{$\Box_{a,b}$} for the rectangular partition with $a$ rows, each of length $b{-}a$.
For example, 
\[
\Box_{1,2}\ =\  \raisebox{-2pt}{\includegraphics{pictures/1}}\,,\quad
\Box_{1,6}\ =\  \raisebox{-2pt}{\includegraphics{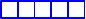}}\,,\quad
\Box_{2,4}\ =\  \raisebox{-6pt}{\includegraphics{pictures/22}}\,,\quad \mbox{ and } \quad 
\Box_{3,7}\ =\  \raisebox{-10pt}{\includegraphics{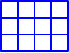}}\ .
\]
Every Schubert problem in Example~\ref{Eq:Derksen-Vakil} has the form $\Box_{a,b}^{\, 4}$.

Fix integers $1\leq a <b$ and $1\leq m<n$, and let $\muDot=(\mu^1,\dotsc,\mu^s)$ be any Schubert problem on $\Gr(m,n)$.
Set $\defcolor{r}\vcentcolon=n{-}m$.
The composition, \defcolor{$\Box_{a,b}^{\, 4}\circ\muDot$},  of $\Box_{a,b}^{\, 4}$ and $\muDot$ is the Schubert problem on $\Gr(2a+m,2b+n)$
given by the following partitions
\[
   \Box_{a,b+r}\,,\, \Box_{a,b+r}\ ,\ 
   \Box_{a+m,b+m}\,,\,  \Box_{a+m,b+m}\ ,\
   \mu^1,\dotsc,\mu^s\,.
\]
Suppose that $a=1$, $b=2$, and $\muDot=\raisebox{-.5pt}{\includegraphics{pictures/1s}}^4$.
  Then $m=r=2$ and $n=4$, so that these partitions are 
\[
   \Box_{1,4}\,,\, \Box_{1,4}\ ,\ 
   \Box_{3,4}\,,\,  \Box_{3,4}\ ,\
   \raisebox{-.5pt}{\includegraphics{pictures/1s}},
   \raisebox{-.5pt}{\includegraphics{pictures/1s}},
   \raisebox{-.5pt}{\includegraphics{pictures/1s}},
   \raisebox{-.5pt}{\includegraphics{pictures/1s}}\,,
\]
which is the Schubert problem $\kapDot$~\eqref{Eq:firstComposition} of Example~\ref{Ex:First_Composition}.

\begin{proposition}[Theorem~3.8 of~\cite{SWY}]
  \label{P:composition}
  The Schubert problem $\Box_{a,b}^{\, 4}\circ\muDot$ is fibered over the Schubert problem $\Box_{a,b}^{\, 4}$ on
  $\Gr(2a,2b)$ with fiber the Schubert problem $\muDot$ on $\Gr(m,n)$.
  We have
  \[
    \Gal_{\Box_{a,b}^{\, 4}\circ\muDot}\ \subset\ \left(\Gal_{\muDot}\right)^{\binom{b}{a}} \rtimes S{\textstyle\binom{[b]}{a}}\,.
  \]
\end{proposition}

We conjecture this inclusion is an equality---that would prove the existence statement in Conjecture~\ref{C:IGP}, for Grassmannians. 

\begin{proof}[Sketch of proof.]
  We explain how to decompose a general instance of the Schubert problem $\Box_{a,b}^{\, 4}\circ\muDot$.
  This is similar to  Example~\ref{Ex:First_Composition}.
  This involves constructing an instance of $\Box_{a,b}^{\, 4}$ as an auxiliary problem, and for each of its solutions, constructing an
  instance of $\muDot$.
   
 An instance of the Schubert problem $\Box_{a,b}^{\, 4}\circ\muDot$ is given by two $b$-planes $K^1_b,K^2_b$, two codimension
 $b$-planes $L^3_{b+n},L^4_{b+n}$, and flags $F^1,\dotsc,F^s$ in $\bbk^{2b+n}$.
 The Schubert problem seeks the $(2a{+}m)$-planes $H$ such that for every $i=1,2$, $j=3,4$, and $t=1,\dotsc,s$, we have 
\[
   \dim H\cap K^i_b\ \geq\ a\,,\ \ 
   \dim H\cap L^j_{b+n}\ \geq\ a+m\,,\ \mbox{ and }\ 
   H\in\Omega_{\mu^t}F^t\,.
\]

 Suppose that the linear subspaces $K^i_b, L^j_{b+n}$, and the flags $F^t$ are in general position.
 Let $H$ be a solution to this instance of $\Box_{a,b}^{\, 4}\circ\muDot$.
 If we set $\defcolor{\Lambda}\vcentcolon=K^1_b\oplus K^2_b\simeq\bbk^{2b}$ and $\defcolor{M}\vcentcolon=L^3_{b+n}\cap L^4_{b+n}\simeq\bbk^n$, then
 $\dim H\cap\Lambda=2a$ and $H\cap M=m$.
 Setting $\defcolor{K^j_b}\vcentcolon=\Lambda\cap L^j_{b+n}$ for $j=3,4$, we have that $\defcolor{h}\vcentcolon=H\cap\Lambda$ is a solution to the Schubert
 problem $\Box_{a,b}^{\,4}$ in $\Gr(2a,\Lambda)$ given by $K^1_b,\dotsc,K^4_b$.
 Let $1\leq t\leq s$.
 As the flag $F^t$ is in general position with respect to the linear spaces $K^1_b$, $K^2_b$, $L^3_{b+n}$, and $L^4_{b+n}$, it is in general
 position with respect to $\Lambda$ and $h$.
 Consequently, $\dim(h+F^t_r)=\dim(h)+\dim(F^t_r)=2a+r$.
 As $M$ has codimension $2b$ and also meets $h+F^t_r$ properly, $(h+F^t_r)\cap M$ has dimension $2a{+}r{-}2b$.
 Thus for $1\leq r\leq n$, $(h+F^t_{r+2b-2a})\cap M$ defines a flag \defcolor{$F^t(h)$} in $M$.
 A further exercise in dimension-counting and the definition of Schubert variety~\eqref{Eq:Schubert_variety} shows that
 $H\cap M\in\Omega_{\mu^t} F^t(h)$.

 Furthermore, for every solution $h$ to the auxiliary problem $\Box_{a,b}^{\,4}$ in $\Gr(2a,\Lambda)$, if we define flags $F^i(h)$ in $M$ 
 as above, then, for every solution $h'$ to the instance of the Schubert problem $\muDot$ in $\Gr(m,M)$ given by the flags
 $F^\bullet(h)$, the direct sum $h\oplus h'$ is a solution to the original Schubert problem.
\end{proof}
  
\subsection{An emerging landscape of Schubert Galois groups}
\label{S:SC_history}

The constructions and results described in Sections~\ref{SS:Enriched_SP} and~\ref{SS:IGP_SP} arose from a sustained investigation of
Schubert Galois groups in which computer experimentation informed theoretical advances.
This began with Vakil's seminal paper~\cite{Va06b}.
There, he used his geometric Littlewood-Richardson rule~\cite{Va06a} in a method that can show a Schubert Galois group is giant.
He applied this method to study Schubert Galois groups in small Grassmannians.
Subsequent experimentation and results this inspired is leading to an understanding what to expect for Schubert Galois groups, and an
outline of a potential classification is emerging from this study.

Without delving into its (considerable) details, we sketch salient features of Vakil's geometric Littlewood-Richardson
rule~\cite{Va06b}.
Given a Schubert problem $\muDot$ on a Grassmannian $G(m,n)$, it constructs a tree \defcolor{$\calT_{\muDot}$}  with $d(\muDot)$ leaves that
encodes a sequence of deformations of intersections of Schubert varieties as the flags move into special position.
Each node $\cbd$ of $\calT_{\muDot}$ determines an enumerative problem which involves intersecting a subset of the Schubert
varieties in $\muDot$ with a \demph{checkerboard variety} \defcolor{$Y_{\cbd}(E,M)$}.
Here $E,M$ are two flags in a special position (determined by $\cbd$) and $Y_{\cbd}(E,M)$ is the set of $m$-planes in $G(m,n)$ having a
particular position with respect to $E,M$ (again specified by $\cbd$).
Let \defcolor{$d(\cbd)$} be the number of solutions to this enumerative problem and  \defcolor{$\Gal(\cbd)$} be its Galois group.

The root of the tree $\calT_{\muDot}$ is labeled by $\muDot$.
For a leaf node $\cbd$, $d(\cbd)=1$.
Every node $\cbd$ of $\calT_{\muDot}$ that is not a leaf has either one child $\cbd'$ or two children $\cbd'$ and $\cbd''$, and we have
$d(\cbd)=d(\cbd')$, respectively $d(\cbd)=d(\cbd')+d(\cbd'')$, when there is one child, respectively two children.
The children of a node are in bijection with the irreducible components of the checkerboard variety $Y_{\cbd}(E,M)$ as the flags $E,M$ become
more degenerate. 

\begin{theorem}[Thms.~3.2 and 3.10 of~\cite{Va06b}]\label{T:Vakil}
  Let $\cbd$ be a node in $\calT_{\muDot}$.
  Suppose that the Galois group of each child of $\cbd$ is
  giant.
  Then $\Gal(\cbd)$ is
  giant
  if one of the following conditions {\rm(1)}, {\rm(2a)}, or {\rm(2b)} hold.
  \begin{enumerate}
    \item $\cbd$ has a unique child.
    \item $\cbd$ has two children $\cbd'$ and $\cbd''$, and
      \begin{enumerate}
        \item  $d(\cbd')\neq d(\cbd'')$ or both are equal to $1$, or
        \item  $\Gal(\cbd)$ is $2$-transitive, and we do not have $d(\cbd')=d(\cbd'')=6$.
      \end{enumerate}
  \end{enumerate}
\end{theorem}

When $\cbd$ is a leaf, $d(\cbd)=1$ so that $\Gal(\cbd)=S_1$ is
giant.
Theorem~\ref{T:Vakil} leads to Vakil's recursive method that may conclude $\Gal(\muDot)$ is
giant.
Given a Schubert problem $\muDot$, this method first constructs $\calT_{\muDot}$, which it then investigates.
If, for every non-leaf node $\cbd$, either (1) or (2a) holds at $\cbd$, then it declares that $\Gal_{\muDot}$ is
giant. 
Otherwise, the method is inconclusive.
It is not a decision procedure, but it is a useful filter to identify Schubert problems that may have enriched Galois groups and thus are
worthy of further study.
The construction of the tree $\calT_{\muDot}$ and the arguments behind Vakil's Theorem~\ref{T:Vakil} hold over any field.

Vakil wrote a Maple script\footnote{\tt http://math.stanford.edu/\~{}vakil/programs/galois} that runs his method on all Schubert problems
on a given Grassmannian. 
He ran this on all small Grassmannians.
Every Schubert Galois group on $\Gr(2,n)$ for $n\leq 16$ and on $\Gr(3,n)$ for $n\leq 9$ was found to be
giant
(for $\Gr(3,n)$, 
Condition (2b) in Theorem~\ref{T:Vakil} was needed).
As $\Gr(4,6)\simeq \Gr(2,6)$ and $\Gr(4,7)\simeq \Gr(3,7)$, the next Grassmannian was $\Gr(4,8)$.
His algorithm was inconclusive for 14 (out of 3501) Schubert problems on $\Gr(4,8)$.
These 14 include the problem $\Box_{2,4}^4$ from Example~\ref{Eq:Derksen-Vakil} and 
the problem of Example~\ref{Ex:First_Composition}.
The Galois groups of these 14 problems are known and none are 2-transitive~\cite[Sect.~5.5]{MSJ}.

A Schubert problem $\muDot$ on $\Gr(m,n)$ is \demph{simple} if at most two of the conditions in $\muDot$ are not
$\raisebox{-.5pt}{\includegraphics{pictures/1s}}=(1,0,\dotsc,0)$.
Using the Pieri homotopy algorithm~\cite{HSS98} to compute solutions to simple Schubert problems and monodromy, Galois groups (over $\CC$)
of many simple Schubert problems (including one with 17,589 solutions) were shown to have full-symmetric Galois groups~\cite{LS09}.
This implies the same for any subfield $\bbk$ of $\CC$.

The first general result concerning Schubert Galois groups was given in~\cite{BMdCS}.
Using Vakil's algorithm and combinatorial reasoning, it was shown that every Schubert problem on $\Gr(2,n)$
for all $n$ has
giant
Galois group.
With an eye towards Condition (2b) in Theorem~\ref{T:Vakil}, another general result showed that Galois groups of Schubert problems on
$\Gr(3,n)$, for every $n$ are 2-transitive~\cite{SW15}.
Liao and Rybnikov showed that all simple Schubert problems with at most one condition
  $\lambda\neq\raisebox{-.5pt}{\includegraphics{pictures/1s}}$, where $\lambda$ is neither symmetric nor a hook, have giant
  Galois group~\cite{LR25}. 
Their method was to study the subfamily of the family $\pi_{\lamDot}\colon \Gamma_{\lamDot}\to \Fl(n)^s$
of~\eqref{Eq:Schubert_branched_cover} given by osculating flags.
When $\lambda$ is symmetric or a hook, these restrictions for osculating Schubert calculus were known~\cite{ModFour,MSJ}.
These results and computations described below suggest the following dichotomy for Schubert Galois groups.

\begin{conjecture}\label{C:dichotomy}
  A Schubert Galois group is either the full symmetric group or it is not 2-transitive.
\end{conjecture}

Robert Williams used the method of computing Frobenius elements to show that many Schubert problems
have full symmetric Galois groups over $\QQ$~\cite{RW}.
These include all Schubert problems on a Grassmannian $\Gr(2,n)$ with up to 500 solutions, as well as all simple Schubert problems on any
Grassmannian with up to 500 solutions, and all Schubert problems on $\Gr(4,9)$ with at most 300 solutions~\cite{GIVIX}.
The numbers here, 300 and 500, are approximate and they represent the limit of the software
used---Singular~\cite{Singular}---to solve a Schubert problem over a prime field (typically $\FF_{1009}$) in a few
hours. 

We close with a sketch of the results from~\cite{GIVIX}, which determined all enriched problems on $\Gr(4,9)$.
This began by using Vakil's method, both his maple implementation and a perl implementation by C.~Brooks, to
identify many Schubert problems on $\Gr(4,9)$
whose Galois group is giant.
For only 233 of the 38,760 Schubert problems was the method inconclusive, and further study found exactly 149 Schubert problems on $\Gr(4,9)$
that had enriched Galois groups.
We remark that this (and earlier computations on $\Gr(4,8)$) only tested Schubert problems that could not
be reduced to a Schubert problem on a smaller Grassmannian. 

Each of these 149 enriched Schubert problems was shown to be a fibration as in Definition~\ref{D:fiberedSchubertProblem},
where the constituent Schubert problems were on a $\Gr(2,4)$ or a $\Gr(2,5)$, and had full symmetric Galois groups, either $S_2$ or $S_3$ or
$S_5$.
By Proposition~\ref{P:wreath}, the Schubert Galois group of each was a subgroup of a wreath product of symmetric groups.
Computing sufficiently many Frobenius elements showed that in each case, the Galois group over $\QQ$ was the expected wreath product.
This computation is explained and archived on the web
page\footnote{{\tt https://FrankSottile.github.io/research/stories/GIVIX/}}.
Of these, 120 are compositions of Schubert problems as in Proposition~\ref{P:composition}, while the remaining 29 had a
different structure. 
All were shown to have the expected Galois groups over $\CC$, but the arguments given in~\cite{GIVIX} hold over any field.

While these results on Schubert Galois groups have not resulted in a classification, there is an emerging landscape of what
to expect, which we summarize for Grassmannians $\Gr(m,n)$.

\begin{itemize}

  \item If $\min\{m,n{-}m\}=1$, then $\Gr(m,n)\simeq\PP^{n-1}$, and Schubert calculus becomes linear algebra; all Schubert
    problems have one solution.
  There are no non-trivial Galois groups.

\item If $\min\{m,n{-}m\}=2$, then  $\Gr(m,n)\simeq \Gr(2,n)$ and all Schubert Galois groups are
  giant~\cite{BMdCS} and are conjectured to be fully symmetric.

  \item If $\min\{m,n{-}m\}=3$, then  $\Gr(m,n)\simeq \Gr(3,n)$ and all Schubert Galois groups are 2-transitive~\cite{SW15}
    and are conjectured to be fully symmetric.

  \item If $\min\{m,n{-}m\}\geq 4$, then $\Gr(m,n)$ has enriched Schubert problems.
        We conjecture that an enriched Schubert problem is either equivalent to one of Vakil's problems from
        Example~\ref{Eq:Derksen-Vakil}, 
        %
        %
        or it is a fibration of Schubert problems.

\end{itemize}

We also remark that while we do not know whether or not Schubert Galois groups depend upon the base field, we
conjecture that they do not.
Also, this study has barely started in arbitrary Lie types, but as mentioned in Remark~\ref{R:LGn}, it is already
  more complex.


\section{Galois groups in applications}\label{Sec:Applications}
Structures in polynomial systems or in enumerative geometry give information about the associated Galois groups.
In a growing number of applications of algebraic geometry, information about associated Galois
groups may be used to detect these structures and then exploit them for solving or for understanding the application.
We sketch this in three application realms.
Esterov's partial determination of Galois groups for sparse polynomial systems leads to a surprisingly efficient algorithm
to recursively decompose and solve sparse systems.
Work in vision reconstruction problems uses Galois groups to detect decompositions, which are then used in efficient
solvers. 
The classical problem of Alt, to determine four-bar mechanisms whose coupler curve passes through nine given points, has a
hidden symmetry of order six coming from the structure of the problem and its formulation as a system of equations.

\subsection{Solving decomposable sparse polynomial systems}
By Proposition~\ref{P:DecomposableIsImprimitive},  when a branched cover $\pi\colon X\to Z$ is decomposable in that there
is a Zariski open subset $V\subset Z$ over which $\pi$ factors, 
 \begin{equation}\label{Eq:newFactor}
   \pi^{-1}(V)\ \xrightarrow{\; \varphi\;}\ Y\ \xrightarrow{\; \psi\;}\ V\,, 
 \end{equation}
then its Galois group $\Gal_\pi$ is imprimitive (and {\it vice-versa}).
Am\'endola, Lindberg, and Rodriguez~\cite{ALR}
proposed methods to exploit this structure in numerical algebraic geometry.
For example, when the decomposition~\eqref{Eq:newFactor} is known explicitly, fibers of $\pi\colon X\to Z$ may be recovered
from the partial data consisting of one fiber of $\varphi\colon\pi^{-1}(V)\to Y$ and one fiber of $\psi\colon Y\to V$.
They illustrated this on examples where $V=Z$ and the first map $\varphi\colon X\to Y$ comes from the invariants of a group
acting on the fibers of $\pi$.
Interestingly, their methods do not require knowledge of the full Galois group, only of the
decomposition~\eqref{Eq:newFactor}.

This approach is particularly fruitful for the sparse polynomial systems of Section~\ref{S:GGsparse}, whose notation and definitions we use.
Let $\calAdot$ be a collection of supports for a sparse polynomial system.
By Esterov's Theorem~\ref{Th:Esterov}, if $\calAdot$ is either strictly lacunary or strictly triangular, then $\Gal_{\calAdot}$ is
imprimitive, and $\Gal_{\calAdot}$ is completely determined (either a group of units or full symmetric) in all other cases.
Not only does this characterize when the branched cover $\pi\colon\Gamma_{\calAdot}\to \CC^{\calAdot}$ is decomposable, but it leads to an
algorithmic procedure for an explicit description of the decomposition.
We sketch that; A complete description is given in~\cite{DSS}.

When $\calAdot$ is lacunary,  $\ZZ^n/\ZZ\calAdot$ is a nontrivial finite group.
Let $\defcolor{\varphi}\colon(\CC^\times)^n\twoheadrightarrow(\CC^\times)^n$ be the map induced by the inclusion
$\ZZ^n\xrightarrow{\sim}\ZZ\calAdot\subset\ZZ^n$
and the functor $\Hom(\bullet,\CC^\times)$.
This has kernel  $\Hom(\ZZ^n/\ZZ\calAdot,\CC^\times)$, a group of units in $(\CC^\times)^n$.
If \defcolor{$\mathcal{B}_\bullet$} is the preimage of $\calAdot$ under the identification
$\ZZ^n\xrightarrow{\sim}\ZZ\calAdot$, then a system $F(x)$ with support $\calAdot$ has the form $G(\varphi(x))$, where the
system $G$ has support $\mathcal{B}_\bullet$. 
Thus $\Hom(\ZZ^n/\ZZ\calAdot,\CC^\times)$ acts on the solutions of any system $F\in \CC^{\calAdot}$, and in fact on the
branched cover $\pi\colon\Gamma_{\calAdot}\to \CC^{\calAdot}$.
This action is free, and the quotient variety is
the branched cover $\psi\colon\Gamma_{\mathcal{B}_\bullet}\to \CC^{\mathcal{B}_\bullet}(=\CC^{\calAdot})$.
Thus we have the factorization
 \begin{equation}\label{Eq:sparseFactor}
   \pi\ \colon\ \Gamma_{\calAdot} \ \xrightarrow{\ \varphi\ }\
   \Gamma_{\mathcal{B}_\bullet}\  \xrightarrow{\ \psi\ }\ \CC^{\mathcal{B}_\bullet}=\CC^{\calAdot}\,.
 \end{equation}
We have $\MV(\calAdot)=|\ZZ^n/\ZZ\calAdot|\cdot\MV(\mathcal{B}_\bullet)$ and when $\MV(\calAdot)>|\ZZ^n/\ZZ\calAdot|$, the
decomposition of $\pi\colon\Gamma_{\calAdot}\to\CC^{\calAdot}$ through the intermediate variety $\Gamma_{\mathcal{B}_\bullet}$ is
nontrivial. 

The first of the maps in this decomposition is induced from the inclusions
$\Gamma_{\calAdot}\subset\CC^{\calAdot}\times(\CC^\times)^n$ and 
$\Gamma_{\mathcal{B}_\bullet}\subset\CC^{\mathcal{B}_\bullet}\times(\CC^\times)^n$ by the identification
$\CC^{\calAdot}=\CC^{\mathcal{B}_\bullet}$ and the monomial map $\varphi\colon(\CC^\times)^n\to(\CC^\times)^n$, from which
fibers may be explicitly computed.  
The map $\varphi$ is computed from the Smith normal form of a matrix whose columns are generators of $\ZZ\calAdot$.
Similarly, the second map is simply the branched cover associated to the family of sparse systems of support $\mathcal{B}_\bullet$.
If $\mathcal{B}_\bullet$ is lacunary or triangular, then this may be further decomposed.
If not, then its fibers are readily computed by numerical software such as \texttt{PHCpack}~\cite{PHCpack} or
\texttt{HomotopyContinuation.jl}~\cite{BT}, using the polyhedral homotopy~\cite{HS95}.

Suppose now that $\calAdot$ is strictly triangular with witness $\emptyset\neq I\subsetneq [n]$, so that $\rank(\ZZ\calA_I)=|I|$
and $1<\MV(\calA_I)<\MV(\calAdot)$.
Then there is a monomial change of coordinates on $(\CC^\times)^n$ and a reindexing of the supports so that $I=[m]$ 
and $\calA_1,\dotsc,\calA_m\subset \ZZ^m$, which is the first $m$ coordinates of $\ZZ^n$.
Writing $(\CC^\times)^n = (\CC^\times)^m\times(\CC^\times)^{n-m}$ for the corresponding splitting, points $x\in(\CC^\times)^n$
are ordered pairs $x=(y,z)$ with $y\in(\CC^\times)^m$ and $z\in (\CC^\times)^{n-m}$.
Then a system $F$ with support $\calAdot$ has the form $F(x)=(G(y),H(y,z))$, where $G$ has support $\calA_I$ and $H$ has support
$\calA_{I^c}$, where $\defcolor{I^c}\vcentcolon=\{m{+}1,\dotsc,n\}$.
Any solution to $F=0$ is a pair $(y^*,z^*)$, where $y^*$ is a solution to $G(y)=0$, and $z^*$ is a solution to the system
$H(y^*,z)=0$ on $(\CC^\times)^{n-m}$.
This structure is apparent in the decomposition
 \begin{equation}\label{Eq:TrianDecomp}
   \Gamma_{\calAdot}\ \longrightarrow\ \Gamma_{\calA_I}\times\CC^{\calA_{I^c}}
   \ \longrightarrow\ \CC^{\calA_I}\times\CC^{\calA_{I^c}}\ =\ \CC^{\calAdot}\,,
 \end{equation}
where the first map is induced by the projection $(\CC^\times)^n\to(\CC^\times)^m$ onto the first $m$ coordinates applied to solutions $(y^*,z^*)$.

Let $p\colon \ZZ^n\to\ZZ^{n-m}$ be the projection to the last $n{-}m$ coordinates and observe that for any solution $y^*$ to $G(y)=0$,
$H(y^*,z)$ has support $p(\calA_{I^c})$.
We have the following product formula (see~\cite[Lem.~6]{ThSt} or~\cite[Thm.~1.10]{Esterov}),
\[
  \MV(\calAdot)\ =\ \MV(\calA_I)\cdot\MV(p(\calA_{I^c}))\,.
\]
Since $1<\MV(\calA_I)$ and $1<\MV(\calAdot)/\MV(\calA_I)=\MV(p(\calA_{I^c}))$, the decomposition~\eqref{Eq:TrianDecomp}
is nontrivial. 
When either $\calA_I$ or $p(\calA_{I^c})$ is lacunary or triangular, these maps may be further decomposed.
If one is neither lacunary or triangular, then its fibers are readily computed by numerical software as above.

\begin{figure}[htb!]
	\begin{center}
	{\includegraphics[scale=.26]{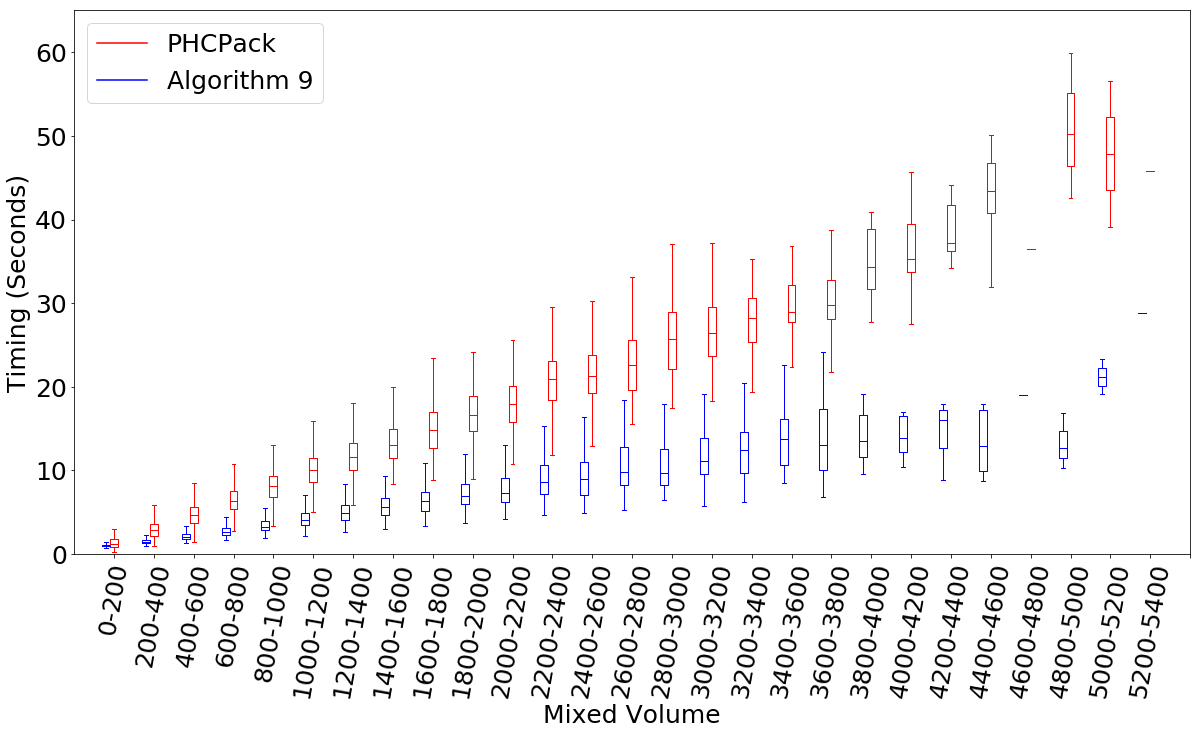}}
   \end{center}
\caption{Box plot of timings comparing \texttt{PHCpack} and Algorithm 9.}
\label{F:Box}
\end{figure}

This leads to an algorithm to recursively decompose the branched
cover $\pi\colon\Gamma_{\calAdot}\to\CC^{\calAdot}$.
In each decomposition, the decomposability of each factor is determined by examining another sparse family.
A blackbox solver (e.g.~\texttt{HomotopyContinuation.jl}~\cite{BT} or \texttt{PHCpack}~\cite{PHCpack}) to compute fibers of
indecomposable branched covers,
combined with the methods of Am\'endola, et.\ al~\cite{ALR}, results in an efficient algorithm for solving sparse polynomial
systems, which was developed in~\cite{DSS}.
These methods have been implemented in the Macaulay2~\cite{M2} package \texttt{DecomposableSparseSystems.m2}~\cite{DSSM2}.
In~\cite{DSS} this package was used in an experiment in which thousands of decomposable systems were solved, both using the
black box solver \texttt{PHCpack} and that package (called Algorithm 9 in~\cite{DSS}).
Figure~\ref{F:Box} shows a box plot of the timings.
This was Fig.~2 in~\cite{DSS}.

\subsection{Vision Problems}
A camera takes a 2-dimensional image of a 3-dimensional scene.
The fundamental problem of image reconstruction is to recover the scene from images taken by cameras at different unknown
locations.  
For this, some features (e.g.\ points, lines, and incidences) are matched between images.
This matching is used to infer the camera positions, which are then used for the full reconstruction.
There are many versions of this nonlinear problem of determining camera positions---different types of cameras and different configurations
of matched features.

A calibrated perspective camera consists of a \demph{focal point} $\bft\in\RR^3$ and a direction vector $\bfv$.
The image is the projection of $\RR^3{\smallsetminus}\{\bft\}$ from the point $\bft$ onto a plane with normal vector $\bfv$ lying a distance 1
from $\bft$ in the direction of $\bfv$.
The image of a point $\bfx\in\RR^3$ is the intersection of the line between $\bft$ and $\bfx$ with this plane.
The considered features are some points, lines, and their incidences, which are assumed to be present in each image.

An \demph{image reconstruction problem} is specified by the number of cameras (images) and the matched features.
For example, we may have two cameras and five points in each image.
Such a problem is \demph{minimal} if, for general data, there is a positive, finite number of solutions (camera positions).
The \demph{degree} of the minimal problem is this number of (complex) solutions for general data, which is a measure of the algebraic
complexity of solving the minimal problem.
Highly optimized solvers have been developed for some minimal problems~\cite{Kukelova,Nister}.
The minimal reconstruction problems were recently classified~\cite{PLMP}, finding many new minimal problems.
Among these new minimal reconstruction problems are some which have imprimitive Galois groups, whose corresponding
decomposable structure may be exploited for solving~\cite{ALR}.

We  present some of the formulation of reconstruction problems.
Fix a reference frame, choosing one camera to be at the origin and to face upwards.
Any other camera is the translation of the first by an element of the special Euclidean group, $\SE_\RR(3)$.
A element of $\SE_\RR(3)$ is a pair \defcolor{$[\bfR\mid\bft]$}, where $\bfR\in\SO(3)$ is a rotation matrix and $\bft\in\RR^3$ is a
translation vector.
Then $[\bfR\mid\bft]$ represents a camera with focal point $\bft$ and direction vector $\bfR\bfk$, where $\bfk$ is the upward-pointing unit
vector. 
In this way, elements of $\SE_\RR(3)$ give coordinates for cameras.
The fixed camera has coordinate $[\bfI\mid\bfzero]$ where \defcolor{$\bfI$} is the identity matrix and \defcolor{$\bfzero$} is the zero
vector. 

\begin{figure}[htb]
\begin{picture}(197,117)(-7,0)
   \put(0,0){\includegraphics{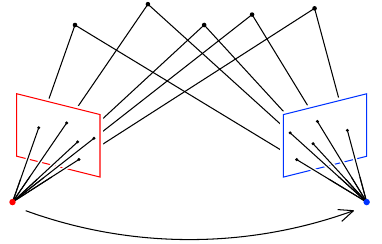}}
   \put(-7,8){$C_1$} \put(177,8){$C_2$}
      \put(77,10){$[\bfR\mid\bft]$}
\end{picture} 
\caption{Minimal problem of two cameras with five points.}
\label{F:twoCameras}
\end{figure}
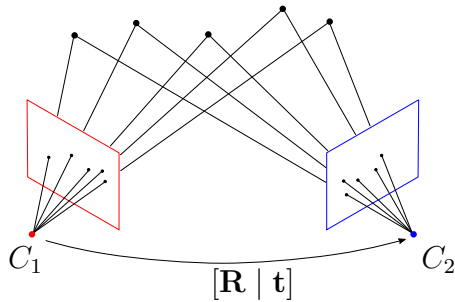

The image plane \defcolor{$\Pi$} of a camera $[\bfR\mid\bft]$ consists of the points $\bfp\in\RR^3$ satisfying the equation
$(\bfR\bfk)\cdot(\bfp-\bft)=1$.
For $\bfx\in\RR^3\smallsetminus\{\bft\}$, its image in $\Pi$ is the point
\[
   \bft\ +\ \frac{\bfx-\bft}{(\bfR\bfk)\cdot(\bfx-\bft)}\ .
\]
Translating by $-\bft$ and applying $\bfR^{-1}$ sends the image plane $\Pi$ to the standard reference plane
$\defcolor{\Pi_{\bfzero}}\vcentcolon=\{(x,y,1)\mid x,y\in\RR\}$ for the
camera $[\bfI\mid\bfzero]$.
We use the coordinates from $\Pi_{\bfzero}$ to represent images of points for all camera.
Thus a point $\bfy\in\Pi_{\bfzero}$ is the image of a point $\bfx\in\RR^3$ under the camera $[\bfR\mid\bft]$ if
 \begin{equation}\label{Eq:Images}
   \bfx\ =\ \bfR \alpha \bfy\ +\ \bft\,,
 \end{equation}
where $\alpha =(\bfR\bfk)\cdot(\bfx-\bft)$ is the \demph{focal depth} of the point $\bfx$ relative to $[\bfR\mid\bft]$.
Figure~\ref{F:twoCameras} is a schematic showing the correspondence between five points $\bfx\in\RR^3$ and their images in the planes $\Pi$,
for two cameras.

Given matched configurations of points, lines, and incidences in $\Pi_{\bfzero}$ for each of several, say $n$, cameras,
equations based 
on~\eqref{Eq:Images} formulate the image reconstruction problem as a system of equations on $(\SE_\RR(3))^{n-1}$.
Complexifying gives a system of polynomials that depends upon the input configuration.
When the problem is minimal, this gives a branched cover over the parameter space of all input configurations.
The degree of the branched cover is the degree of the minimal problem.
As we have seen before, there is a Galois group for each minimal problem.
When the Galois group is imprimitive, Proposition~\ref{P:DecomposableIsImprimitive} implies that the branched cover is
decomposable. 
If a decomposition~\eqref{Eq:newFactor} is known, then that may be exploited for solving.


One such problem with imprimitive Galois group is that of reconstructing five points given images from two cameras, which is illustrated in
Figure~\ref{F:twoCameras}. 
The branched cover corresponding to this minimal problem has degree 20.
The imprimitivity may be understood by observing that the solutions come in pairs:
Given one solution $([\bfI\mid\bfzero],[\bfR\mid\bft])$, a second is given by rotating the camera $[\bfR\mid\bft]$ $180^\circ$
around the line between the two cameras.
(This also changes the inferred positions of the unknown points $\bfx\in\RR^3$.)
This is called a \demph{twisted pair} in the literature, and we see that the Galois group preserves the resulting partition of the 20
solutions into ten twisted pairs, and is hence a subgroup of $S_2\wr S_{10}$.
In fact, the Galois group is even smaller, it is $(S_2\wr S_{10})\cap A_{20}$~\cite{VisMono}, which is the Weyl group $D_{10}$.
This imprimitivity implies the associated branched cover is decomposable and the system can be solved in stages.
A decomposition for this problem is implicit in~\cite{Nister}.

In~\cite{VisMono}, the minimal problems of degree at most 1000 with imprimitive Galois group were classified.
Those were further studied using numerical algebraic geometry, which led to an understanding of their structure, and for
many an explicit decomposition was found.

\subsection{Alt's Problem}
Polynomial systems arise in engineering when designing mechanisms with a desired range of motion.
Robotic arm movements, for instance, may need
to reach several positions to perform specific tasks.
These movements can be modeled by polynomial systems, from which they can be studied with the methods discussed.
One such problem due to Alt~\cite{Alt} is the nine-point synthesis problem for four-bar linkages.

A \demph{four-bar linkage} is a planar mechanism built from a quadrilateral (which may self-intersect) with rotating joints and fixed
side lengths.
One side of the quadrilateral is the \demph{base} which is fixed in place, while the other sides move as allowed by freely rotating
the joints.
The side opposite the base is the \demph{coupler bar}, and a triangle is erected on the coupler bar.
In an actual mechanism, a tool is placed at the apex of the triangle and the mechanism is
maneuvered to position the tool.

\begin{figure}[htb]
  \centering
  \includegraphics[height=120pt]{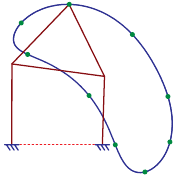}
  \caption{A four-bar mechanism and coupler curve.}
  \label{Fig:four-bar}
\end{figure}

To understand this motion, consider the quadrilateral.
Removing the coupler bar, the two bars that were incident to it may each rotate freely around their fixed points.
The coupler bar imposes a distance constraint on the rotating bars, and there remains one degree of freedom.
(The abstract curve of this motion has genus one.)
In the resulting motion, the apex of the triangle traces the \demph{coupler curve}, see Figure~\ref{Fig:four-bar}.

The space of all mechanisms is nine-dimensional.
Indeed, the two fixed points may be any points in $\RR^2$, giving four dimensions (degrees of freedom).
The lengths of each of the remaining five segments in the mechanism give five more, for a total of nine.
That the coupler curve contains a given point in $\RR^2$ is a single, simple condition on the space of four-bar mechanisms.
Thus we expect there are only finitely many mechanisms whose coupler curve contains nine given general points.
Alt~\cite{Alt} recognized this, and his nine-point synthesis problem asks for the mechanisms whose coupler curve contains
a given nine points.

Identifying $\RR^2$ with $\CC$, we represent the bars as complex numbers.
Complexifying gives a useful formulation of Alt's problem in \demph{isotropic coordinates}---this is described
in~\cite{MSW}. 
Solutions to these equations for nine general points were computed using homotopy continuation in~\cite{MSW}, finding 8652
solutions. 

In~\cite{HS12}, this computation was repeated and a soft certificate was computed to certify
the 8652 computed solutions.
While 8652 is almost surely the number of solutions, these computations only show that it is a lower bound, and a proof of
the number 8652 remains elusive.
Further evidence for the number 8652 was found in~\cite{ProbSat}, but that result also only implies that 8652 is a
lower bound.

In this formulation, solutions come in pairs due to relabeling---swapping labels of the bars incident to the base and to
the apex of the triangle results in another solution and gives the same four-bar mechanism.
Classically Roberts~\cite{Roberts} and Chebyshev~\cite{MaSo} (see~\cite{Ver} for a discussion) show that there are three
mechanisms---called Robert's Cognates---with the same coupler curve.
Consequently, the Galois group of this formulation of Alt's problem is imprimitive as it preserves the six solutions which
give the same coupler curve.
We also see that, assuming the number 8652 is correct, that there are 4326 four-bar mechanisms whose coupler curve contains
nine given  points, and 1442 distinct coupler curves.

Since label swapping may be done independently on each cognate, the Galois group \defcolor{$G$} of the six solutions with
given coupler curve is $\ZZ/2\ZZ\times\ZZ/3\ZZ=\ZZ/6\ZZ$.
Consequently, the Galois group of this formulation is a subgroup of $(\ZZ/6\ZZ)\wr S_{1432}$.
To the best of our knowledge, the Galois group of this problem has not been determined.

\bibliographystyle{amsplain} 
\providecommand{\bysame}{\leavevmode\hbox to3em{\hrulefill}\thinspace}
\providecommand{\MR}{\relax\ifhmode\unskip\space\fi MR }
\providecommand{\MRhref}[2]{%
  \href{http://www.ams.org/mathscinet-getitem?mr=#1}{#2}
}
\providecommand{\href}[2]{#2}

\end{document}